\documentclass[smallextended,envcountsect]{svjour3} 

\usepackage{graphicx,algorithm,color}
\usepackage[legalpaper, margin=1in]{geometry}
\usepackage{amsmath}
\usepackage{framed}
\usepackage{amssymb}
\usepackage{subfigure,verbatim}
\newtheorem{assumption}{Assumption}[section] 
\newcounter{algo}[section]

\def\xkd{x_k^{\delta}}
\def\xk+1d{x_{k+1}^{\delta}}

\newcommand{\ja}{J(\xkd)^T}

\newcommand{\Jdk}{J(\xkd)}

\newcommand{\TT}{T_{\ell_k}^TT_{\ell_k}}
\newcommand{\ql}{Q_{\ell_k}}
\newcommand{\tl}{T_{\ell_k}}
\newcommand{\A}{T_{\ell_k}^TT_{\ell_k}}
\newcommand{\ke}{k_\epsilon}

\begin{document}

\title{An inexact non stationary Tikhonov procedure for large-scale nonlinear ill-posed problems }
\author{S. Bellavia, M.Donatelli, E.Riccietti}

\author{Stefania  Bellavia  \and Marco Donatelli \and Elisa Riccietti 
}


\institute{S. Bellavia \at Dipartimento  di Ingegneria Industriale, Universit\`a di Firenze,
	viale G.B. Morgagni 40,  50134 Firenze,  Italy. 
             \and
M. Donatelli \at  	Dipartimento di Scienza e Alta Tecnologia, Universit\`a dell' Insubria, 
via Valleggio 11, 22100 Como, Italy.	
           \and
           E. Riccietti \at
              IRIT-INPT, rue Charles Camichel, B.P. 7122 31071, Toulouse Cedex 7, France.
       \email{elisa.riccietti@enseeiht.it} 
}

\date{Received: date / Accepted: date}


\maketitle

\begin{abstract}
In this work we consider the stable numerical solution of  large-scale ill-posed nonlinear least squares problems with nonzero residual. We propose a non-stationary Tikhonov  method with inexact step computation, specially designed  for large-scale problems. At each iteration the method requires the solution of an elliptical trust-region subproblem to compute the step. This task is carried out employing a Lanczos approach, by which an approximated solution is computed. The trust region radius is chosen to ensure the resulting Tikhonov regularization parameter to  satisfy a prescribed condition on the model, which is proved to ensure regularizing properties to the method. The proposed approach is tested on a parameter identification problem and on an image registration problem, and it is shown to provide important computational savings with respect to its exact counterpart.
\end{abstract}

\section{Introduction}
Let $F:\mathcal{D}(F)\subseteq\mathbb{R}^n\rightarrow\mathbb{R}^m$, with $m\geq n$, be a continuously differentiable nonlinear function.  We denote with $\mathcal{D}(F)$ the domain of definition of $F$  and consider problems of the form: 
\begin{equation}\label{problem}
\min_{x} f(x)=\frac{1}{2}\|F(x)-y\|^2,
\end{equation}
that are ill-posed, in the sense that the solutions do not depend continuously on the data. We assume that a solution  $x^\dagger$ for \eqref{problem} exists. If for all the solutions of \eqref{problem} it holds $\|F(x)-y\|>0$, we say that the problem has nonzero residual.

We suppose to have only noisy data $y^\delta$ at disposal, such that, given $\delta\geq 0$: 
\begin{equation}\label{yd}
\|y-y^\delta\|\leq\delta,
\end{equation}
so that we have to deal with the following noisy problem:
\begin{equation}\label{noisy_pb}
\min_{x} f_\delta(x)=\frac{1}{2}\|F(x)-y^\delta\|^2.
\end{equation}

Many approaches have been proposed in the literature to deal with this problem,  such as nonstationary iterated Tikhonov regularization or regularized versions of trust region methods \cite{nostro,ellittica,hanke,wy1,wy,hanke_CG,pet_tr,weese1993regularization}.

In this paper we assume the problem to be large-scale, i.e. we assume $n$ to be large.  
The large-scale setting has been deeply analysed especially in the case of the linear counterpart of this problem, see for example \cite{donatelli,reichel2007greedy,levin2017stopping,reichel2009error,haber2008numerical,OLS,rojas2002trust,hanke2001lanczos,donatelli2012square,chung2008weighted}, among the others. 
The main approach used in this case combines a standard regularization technique, such as Tikhonov regularization, with an iterative method based on partial Lanczos bidiagonalization of the operator \cite{OLS}. These are called \textit{hybrid methods}  because they apply Tikhonov regularization to the projected problem into the Krylov subspace at each iteration \cite{CNO}. Thus, regularization in hybrid methods is achieved both by Krylov filtering and by appropriate choice of a regularization parameter at each iteration.

In this paper we propose a non-stationary Tikhonov procedure with inexact step computation, specially designed  for nonlinear large-scale ill-posed problems. The method is based on an automatic strategy for choosing the free regularization parameter, that requires the solution of an elliptical trust-region subproblem. The radius of the trust-region is chosen to ensure the resulting Tikhonov regularization parameter to satisfy a prescribed condition on the model, which is proved to ensure regularizing properties to the method. The step to update the iterate is then the solution of the trust-region subproblem, which is approximated employing a Lanczos approach. Hence, our method falls into the class of hybrid methods since at each iteration the linearized problem associated to the Jacobian matrix is solved in the Krylov subspace generated by the partial Lanczos bidiagonalization. 

This method represents an inexact version of the approach proposed in \cite{ellittica}, which is tailored for small-medium scale problems. The method indeed employs the singular value decomposition (SVD)  of the Jacobian matrix, which can be too costly to compute for large-scale problems. Avoiding this computation we introduce several sources of inexactness, that will be described in details in the following section. Due to this, from a theoretical point of view,  we cannot expect the same regularizing properties of its exact counterpart in \cite{ellittica}, but the proposed approach reduces to the method in \cite{ellittica} when the dimension of the Lancsoz space $\ell$ reaches $n$. However, we are still able to show that it is possible to control how the distance  of the current solution approximation from an exact solution of the problem changes from an iteration to the other, and that this decreases if the dimension of the Lanczos space is large enough. 

 In case of  noise-free problems, the proposed procedure is  shown to be globally convergent, as the reduction in the quadratic model provided by the computed step is greater than that provided by the so-called Cauchy step. 
In the case of noisy data, the process is stopped before the convergence is reached, to avoid approaching a solution of the noisy problem. In this case we are not considering the global convergence issue. However, the trust-region mechanism guarantees that the stopping criterion is satisfied after $O(\delta^{-2})$ iterations.

Numerical tests enlighten that in practice the method behaves really well and shows good regularizing properties, just as its exact counterpart. As expected, the proposed approach provides  a considerable overall reduction of the CPU time over its exact counterpart.
Even if theoretical results requires $\ell$ large enough to prove regularization properties, the numerical tests show that such properties are maintained also for small values of $\ell$.

The paper is organized as follows. 
In Section \ref{sec_prelim} we briefly describe the exact method in \cite{ellittica} and the standard Lanczos bidiagonalization method.
In Section \ref{sec_method} we describe the new approach we propose, which represents an extension of the method in \cite{ellittica}, suitable to handle large-scale problems thanks to the use of Lanczos technique in the trust-region scheme.  In Section \ref{sec_radius} we introduce the so-called \textit{projected $q$-condition}, which will be used to prove regularizing properties of the method. We show how to choose the trust-region radius to automatically obtain a regularization parameter that satisfies such condition. Section \ref{sec_algo} and Section \ref{sec_analysis} are devoted to the theoretical analysis of the proposed method, in particular in Section \ref{sec_algo} we discuss the convergence, complexity and computational cost, while in Section \ref{sec_analysis} we focus on the regularizing properties of the proposed method.  Finally in Section \ref{sec_num}, in order to avoid the estimation of the further parameter $\ell$,  we propose to use a nonstationary non-decreasing sequence of values for $\ell$ and 
show the regularizing behavior of the method in practice.

\paragraph{Notations}
Throughout the paper the symbol $\|\cdot\|$ will be used to denote the 2-norm.
We denote the iterates as $\xkd$, if the data are exact $x_k$ may be used in alternative to $\xkd$. 
By $x_0^\delta=x_0$ we denote an initial guess which may incorporate a-priori knowledge of an exact solution.
We also define
\begin{align}
&B_k=J(\xkd)^TJ(\xkd), \quad  g_k=\nabla f_\delta(\xkd)=\ja(F(\xkd)-y^\delta),\label{notation}\\
&m_k(p(\lambda))=B_kp(\lambda)+g_k,\label{def_mod}
\end{align} 
where $J(\xkd)$ is the Jacobian of $F(\xkd)$.   We will denote $\varsigma^k_1\geq\dots\geq\varsigma^k_n$ the singular values of $J(\xkd)$ and $e_1$ the first vector of the canonical basis.

\section{Preliminaries}\label{sec_prelim}
\subsection{Exact method in \cite{ellittica}}
In this section we briefly outline the procedure in \cite{ellittica}. For simplicity we consider the case in which $J$ has full rank, but the procedure can be easily extended also the rank-deficient case, cf. \cite{ellittica}. 

At $k$-th iteration of the method in \cite{ellittica}, given the trust-region radius $\Delta_k>0$ and the current iterate $x_k^{\delta}$,  the following elliptic trust-region subproblem is solved (\cite[ch. 4]{nw}, \cite[\S 7.4]{cgt}):
\begin{equation}\label{TR}
\begin{array}{l}
\displaystyle \min_p \frac{1}{2}\|F(x_k^\delta)-y^\delta+\Jdk p\|^2,\\
\mbox { s.t. } \|B_k^{-1/2} p\|\le \Delta_k ,
\end{array}
\end{equation}
which has a unique solution. 
Letting $z=B_k^{-1/2} p$  problem \eqref{TR} reduces to
\begin{equation}\label{TRz}
\begin{array}{l}
\displaystyle \min_{z}  \frac{1}{2}z^TB_k^2 z+z^TB_k^{1/2}g_k+ f_\delta(x_k^\delta)\\
\mbox { s.t. } \|z\|\le \Delta_k.
\end{array}
\end{equation}
KKT conditions for problem \eqref{TRz} are given by 
\begin{subequations}\label{KKT_z}
	\begin{align}
	&(B_k^2+\lambda I) z(\lambda)=-B_k^{1/2}g_k,\label{eq_passo_z}\\
	&\lambda (\|z(\lambda)\|-\Delta_k)=0\label{system_1},\\
	&\lambda\ge 0,\\
	&\|z(\lambda)\|\le \Delta_k\label{system_end}.
	\end{align}
\end{subequations}
Let $(\lambda_k,z(\lambda_k))\in\mathbb{R}^+\times\mathbb{R}^n$ be the solution of \eqref{KKT_z}, with $z_k=z(\lambda_k)$  minimum norm solution of \eqref{eq_passo_z}. If we let
\begin{equation}\label{plambda}
p(\lambda)=B_k^{1/2}z(\lambda),
\end{equation}
and $p_k=p(\lambda_k)$,  the couple  $(\lambda_k,p_k)$ is  a KKT point for \eqref{TR}.  The solution of \eqref{TR} can be then found by solving \eqref{TRz},  and  through relation \eqref{plambda}.  
At each iteration the trust-region radius is chosen to satisfy:
	\begin{equation}\label{bound_delta}
	\Delta_k\leq\frac{1-q}{\|B_k\|^2}\|B_k^{1/2}g_k\|,
	\end{equation} 
for a given $q\in(0,1)$. The process is stopped at iteration $k^*(\delta)$, when the following discrepancy principle is satisfied:
\begin{equation}\label{discrepance}
\|g_{k^*(\delta)}\|\le \tau_k\delta <\|g_k\|,  \;\;\;0\le k< k^*(\delta),\quad \tau_k=\bar{\tau}\|J(x_k^\delta)\|,
\end{equation}
with $\bar{\tau}$ a chosen constant.

The $k$-th iteration of the process is sketched in Algorithm \ref{algoSVD}.
\begin{algorithm}
\caption{$k$-th iteration of the regularizing trust-region method in \cite{ellittica}}
\label{algoSVD}
		Given $\xkd$,    $\eta\in (0,1)$, $\gamma\in (0,1)$, $0<C_{\min}<C_{\max}$, $q\in (0,1)$. 
		\vskip 1pt
		Exact data:  $y$; Noisy data:  $y^{\delta}$.
		\begin{description}
			\item[1.]  Choose $\displaystyle{\Delta_k}\in \left.\left [C_{\min}\|B_k^{1/2} g_k\|, \, 
			\min\left\{ C_{\max},  \frac{1-q}{\|B_k\|^2}  \right\}\|B_k^{1/2} g_k\|\right]\right.$.
			\item[2.]  Repeat\\ 
			2.1 Find the pair $(z_k,\lambda_k)$ solution of KKT conditions \eqref{KKT_z}.
			\\
			2.2 Set $p_k=B_k^{1/2} z_k$.\\
			2.3  Compute  $${\pi_k(p_k)}=\frac{f_\delta(\xkd)-f_\delta(\xkd+p_k)}{f_\delta(\xkd)-\frac{1}{2}\|F(\xkd)-y^\delta+\Jdk p_k\|^2}.$$ 
			2.4 If $\pi_k(p_k)< \eta$, set  ${\Delta_{k}=\gamma \Delta_k}$.\\
			\hspace*{-20pt}		Until $\pi_k(p_k)\ge \eta$.
			\item [3.] Set  ${x_{k+1}^\delta=\xkd+p_k}$.
		\end{description}
		\end{algorithm}

The acceptance criterion in Step 2.4  provides global convergence of the procedure in case of 
noise-free problems 
 \cite{nw}.  In the case of noisy data, the process is stopped before the convergence is reached, to avoid approaching a solution of the noisy problem. In this case we are not considering the global convergence issue. However,  the trust-region acceptance mechanism guarantees that the stopping criterion is satisfied 
after $O(\delta^{-2})$ iterations.  

Each iteration of the method requires two computations that are expensive if the problem is large scale:
\begin{itemize}
	\item the computation of the square root of the matrix $B_k$, to evaluate the right hand side in \eqref{eq_passo_z}, to compute the step $p(\lambda_k)$ from $z(\lambda_k)$ by \eqref{plambda}, and to update the trust region radius in \eqref{bound_delta},
	\item   the  computation of an approximation to $z(\lambda_k)$  through \eqref{KKT_z} requires the solution of the secular equation (cf. \eqref{secular_eq} below)  via Newton's method \cite{nw,cgt}. This requires the solution of a sequence of linear system of the form \eqref{eq_passo_z}, which are large scale systems.
	\end{itemize}
In \cite{ellittica} these operations are performed computing the SVD decomposition of matrix $J(\xkd)$, which is not feasible if the size of the matrix becomes large. 
Having in mind to extend the procedure to large-scale problems, these computations cannot be performed exactly. It is then necessary to design some cheap procedure to approximate the action of matrix $B_k^{1/2}$ onto a vector, and to find an approximate solution to \eqref{TRz}. For the first task, we will need to define an operator ${\sf sq}:\mathbb{R}^{n\times n}\times\mathbb{R}^n\rightarrow\mathbb{R}^n$, that maps a matrix $A$ and a vector $b$ to an approximation of the product $A^{1/2}b$, so that ${\sf sq}(B_k,g_k)$ will be an approximation to $B_k^{1/2}g_k$. 
Then, the method will produce a sequence of solution approximations $\{\xkd\}$, forming at each iteration the new approximation as $x_{k+1}^\delta=\xkd+p_k$, where $p_k$ is defined as 
\begin{equation}\label{plambda_inexact}
p_k=p(\lambda_k)={\sf sq}(B_k,z(\lambda_k)),
\end{equation} 
with $z(\lambda_k)$ approximated solution of 
\begin{equation}\label{system_appr}
\begin{array}{l}
\displaystyle \min_{z}  \frac{1}{2}z^TB_k^2 z+z^T{\sf sq}(B_k,g_k)+f_\delta(\xkd)\\
\mbox { s.t. } \|z\|\le \Delta_k.
\end{array}
\end{equation}

To avoid the computation of the SVD of $J(\xkd)$ and to build an approximation ${\sf sq}(B_k,z(\lambda_k))$ to $p_k$ and approximately solve the minimization problems \eqref{system_appr}, we employ  the Lanczos method, that we describe in the following section. There are several methods to approximate the SVD of a matrix. We have chosen the Lanczos method as it allows us to cope with both our issues at the same time and because the structure induced in the computed step by the Lanczos method allows us to maintain (and theoretically prove) some of the regularizing properties of the method in \cite{ellittica}, cf. Remark \ref{remark_rank1}.  

\subsection{Lanczos bidiagonalization method}
Lanczos methods are usually used to solve large scale linear systems,  cf. for example \cite{golub1965,gvl}. They make computing the SVD of the operator feasible by projecting the problem onto a subspace of small dimension.  

In practice, given a matrix $A\in\mathbb{R}^{m\times n}$  (generally large) and a scalar $\ell>0$ (generally $l\ll n$), the Lanczos bidiagonalization technique computes 
a sequence of Lanczos vectors $p_j\in \mathbb{R}^m$ and $q_j\in \mathbb{R}^n$ and scalars $\alpha _j$, $\beta_j$, $j=1,2,\ldots,\ell$, by the 
 recursive procedure
  sketched in Algorithm \ref{algo_lanczos}, that is initialized by a vector $q_1\in \mathbb{R}^n$.
 
 \begin{algorithm}
	\caption{GKLB($A$,$q_1$,$\ell$) (Golub-Kahan-Lanczos Bidiagonalization)}
	\label{algo_lanczos}
	Given $A\in\mathbb{R}^{m\times n}$, $q_1\in\mathbb{R}^n$,  $\ell\in\{1,\dots,n\}$. 
	\begin{description}
		\item[1.]  Set $q_1=\frac{q_1}{\|q_1\|}$, $\beta_0=0$.
		\item[2.]\label{ciclo_for}  For $j=1,\dots,\ell$
		\begin{description}
			\item[2.1] $p_j=Aq_j-\beta_{j-1}p_{j-1}$
			\item [2.2] $\alpha_j=\|p_j\|$
			\item [2.3]  $p_j=\frac{p_j}{\alpha_j}$
			\item [2.4]\label{def_q} $q_{j+1}=A^Tp_j-\alpha_jq_j$
			\item[2.5] $\beta_j=\|q_{j+1}\|$
			\item[2.6] $q_{j+1}=\frac{q_{j+1}}{\beta_j}$
		\end{description}
	\end{description}
	\end{algorithm}	

Assume for sake of simplicity that the algorithm is not prematurely halted as  $\alpha_j=0$ or $\beta_j=0$ is met (if $\alpha_j$  ($\beta_j$) equals zero  ones must choose a new vector $p_{j+1}$ ($q_{j+1}$) which is orthogonal to the previous $p_j$'s ($q_j$'s) \cite{golub1965}.)
  Then 
after $\ell$ steps, GKLB has generated the tridiagonal matrix  $T_\ell\in \mathbb{R}^{\ell\times \ell}$ 
$$
T_\ell=\begin{bmatrix}
\alpha_1& \beta_1 &  && \\
 & \alpha_2 & \beta_2  & &\\
 &  & \ddots & \ddots &&\\
 &  &   &\alpha_{\ell-1}&\beta_{\ell-1}\\
 &  &  &  &\alpha_\ell\\
\end{bmatrix},
$$
the  matrix $P_\ell\in \mathbb{R}^{m\times \ell}$, whose $j$-th column is given by vector $p_j$ and the  matrix $Q_\ell \in \mathbb{R}^{n\times \ell}$
with orthonormal columns $q_j$, for $j=1,\ldots,\ell$. 

The first column of $Q_\ell$ is the vector   $\frac{q_1}{\|q_1\|}$ and 
at each step $j$, $j=1,\dots,\ell$, it holds:
\begin{align}
AQ_j&=P_j T_j,\label{jac1}\\
A^TP_j&=Q_j Tj^T+\beta_jq_{j+1}e_{j}^T\label{jac2},
\end{align} 
for $Q_j\in\mathbb{R}^{n\times j}$, $P_j\in\mathbb{R}^{m\times j}$, $T_j\in\mathbb{R}^{j\times j}$, $q_{j+1}\in\mathbb{R}^{n}$, $e_j\in\mathbb{R}^{j}$ is the last column of the identity matrix of size $j$.
It also holds
\begin{equation}\label{cond}
Q_j^TQ_j=I,\quad P_j^TP_j=I, \quad Q_j^Tq_{j+1}=0.
\end{equation}

Algorithm \ref{algo_lanczos} computes the same information as the Lanczos tridiagonalization algorithm applied to the Hermitian matrix $A^TA$. 
Indeed, from \eqref{jac1} and \eqref{jac2} we deduce that for each $j=1,\dots,\ell$:
\begin{align}\label{due_rel}
A^TAQ_j=Q_jT_j^TT_j+\beta_jq_{j+1}e_j^TT_j, &&
Q_j^TA^TA=T_j^TT_jQ_j^T+\beta_jT_j^Te_jq_{j+1}^T.
\end{align}
From these two relations, and taking into account \eqref{cond}, we conclude that:
\begin{equation}\label{diag0}
Q_j^T(A^TA)Q_j=T_j^TT_j.
\end{equation}

In particular, the Lanczos vectors $q_j$ constitute an orthonormal basis of the following Krylov subspace: 
$$
\mathcal{K}_{\ell}(A^TA,q_1)={\sf span}\{q_1,(A^TA)q_1,\dots,(A^TA)^{\ell-1}q_1\}.
$$


If applied to $J(\xkd)$, these approaches are convenient for our problem, as they allow us to solve both our issues.  The factorisation of the form  \eqref{diag0}, provided by the Lanczos process, can be used to approximately evaluate $B_k^{1/2}g_k$, and an approximation to a solution of \eqref{system_appr} can be sought in the Krylov subspace generated by the Lanczos process. We will describe in next section how to combine the two methods described in this section to obtain a hybrid Lanczos method for large-scale ill-posed problems.

\section{Regularizing Lanczos trust-region hybrid approach}\label{sec_method}
At each iteration $k$ of the trust-region method, we apply the Lanczos bidiagonalization technique to matrix $J(\xkd)\in\mathbb{R}^{m\times n}$. 

If we perform $n$ steps of the GKLB procedure, i.e. $\ell=n$,  the Lanzos process produces the following decomposition:
\begin{equation}\label{facto}
J(\xkd)=P_n T_n Q_n^T, \quad B_k=Q_nT_n^TT_n Q_n
\end{equation}
where $P_n$ is a $m\times n$ matrix with orthonormal columns, $Q_n$ is a unitary matrix of order $n$, and $T_n\in \mathbb{R}^{n\times n}$ is an upper bidiagonal matrix.

Employing the factorization \eqref{facto}, given a function $f$, we can evaluate the action of $f(B_k)$ onto a vector $b$ in the following way, cf. \cite[ch.13]{higham2008functions}:
$$
f(B_k)b=f(Q_nT_n^TT_nQ_n^T)b=Q_n f(T_n^TT_n) Q_n^Tb.
$$ 
To compute $f(T_n^TT_n)$ we can take advantage of the special structure of matrix $T_n^TT_n$, that is a tridiagonal matrix. If $q_1=b$ is chosen,  $Q_n^Tb=\|b\|e_1$, where $e_1\in\mathbb{R}^n$ is the first vector of the canonical basis, as a consequence of the orthogonality of matrix $Q_n$.

In particular, if $f(x)=\sqrt{x}$, we can compute the action of $B_k^{1/2}$ onto a vector $b$ as 
$$
B_k^{1/2}b=Q_n (T_n^TT_n)^{1/2} Q_n^Tb.
$$  

However, performing $n$ steps of the Lanczos method and computing the SVD of $T_n$ may be not feasible.  Then, we  perform  $\ell_k<n$ steps of the Lanczos process, obtaining an approximate factorization of matrix $B_k$, as in \eqref{diag0} with $j=\ell_k$ and $A=J(\xkd)$:
\begin{equation}\label{diag1}
Q_{\ell_k}^TB_kQ_{\ell_k}=T_{\ell_k}^TT_{\ell_k}.
\end{equation}
 This can be used to define an approximation to $B_k^{1/2}b$:
 \begin{equation}\label{sqrt}
{\sf sq}(B_k,b):=\ql (\tl^T\tl)^{1/2} \ql^Tb, \quad \forall b\in\mathbb{R}^n.
\end{equation} 
We are actually evaluating the function on the Krylov subspace $\mathcal{K}_{\ell_k}:=\mathcal{K}_{\ell_k}(B_k,q_1)$ and expanding the result  back onto the original space $\mathbb{R}^n$, cf. \cite{higham2008functions}. Of course this is an approximation, that would be exact if $\ell_k=n$. 

Usually a good approximation can be obtained with a $\ell_k\ll n$, then 
$\tl^T\tl$ will have much smaller dimension than $B_k$, and  the evaluation of the function will be feasible by direct computation of the SVD of $\tl$. Then, if we choose $q_1=g_k$, we get the approximation we were looking for:
\begin{equation}\label{tsk}
\tilde{s}_k:={\sf sq}(B_k,g_k)= \ql (\tl^T\tl)^{1/2} \ql^Tg_k.
\end{equation}

The computed SVD can be employed also to find an approximate solution to \eqref{system_appr}, by projecting the problem onto the Krylov subspace generated by the Lanczos factorization. At each iteration $k$ of the trust-region method, we look for a solution in the subspace: $z=\ql w$. The projection of problem \eqref{system_appr} is given by:	\begin{equation*}
	\min_{z\in\mathcal{K}_{\ell_k}} \frac{1}{2}z^TB_k^2z+\tilde{s}_k^Tz+f_\delta(\xkd)\\
           \mbox { s.t. } \|z\|\le \Delta_k.
       \end{equation*}
Noting that    $\|z\|=\|Q_{\ell_k} w\|=\|w\|$, as the columns of $Q_{\ell_k}$ are orthogonal, it    can be reformulated as 
       \begin{equation}\label{sub_w}
	\min_{w\in\mathbb{R}^{\ell_k}} \frac{1}{2}w^TQ_{\ell_k}^TB_k^2Q_{\ell_k}w+w^TQ_{\ell_k}^T\tilde{s}_k+f_\delta(\xkd)\\
           \mbox { s.t. } \|w\|\le \Delta_k.
       \end{equation}
       We remark that from \eqref{cond} and \eqref{due_rel} with $A=J(\xkd)$ and $j=\ell_k$, it follows :
\begin{equation}\label{rank1}
       Q_{\ell_k}^TB_k^2 Q_{\ell_k}=(T_{\ell_k}^TT_{\ell_k})^2+\|q_{{\ell_k}+1}\|^2\beta_{\ell_k}^2T_{\ell_k}^Te_{\ell_k}e_{\ell_k}^TT_{\ell_k}.
\end{equation}
To simplify the subsequent convergence analysis we make a further approximation and we drop the rank-one perturbation, approximating therefore $Q_{\ell_k}^TB_k^2Q_{\ell_k}$ with $(T_{\ell_k}^TT_{\ell_k})^2$, band matrix with bandwidth 2. 
  We highlight that the rank-one term that we choose to neglect is of the form $\|q_{\ell_k+1}\|^2\beta_{\ell_k}^2\alpha_{\ell_k}^2 e_{\ell_k}e_{\ell_k}^T$, and it is zero when $\ell_k=n$. It is therefore a residual term, coming from the fact that we prematurely stop the Lanczos process. As we point out in Remark \ref{remark_rank1}, the presence of this term alters the peculiar structure of the step, which is crucial to prove a fundamental relation between the model and the step itself.


The projected problem becomes then 
\begin{equation}\label{min_proj}
	\min_{w\in\mathbb{R}^{\ell_k}} \frac{1}{2}w^T (T_{\ell_k}^T T_{\ell_k})^2 w+w^TQ_{\ell_k}^T\tilde{s}_k+f_\delta(\xkd)\\
	   \mbox { s.t. } \|w\|\le \Delta_k.
\end{equation}
Reminding that  $Q_{\ell_k}^Tg_k=\|g_k\|e_1$ as we set $q_1=g_k$, 
the KKT  conditions for \eqref{min_proj}  are the following:
\begin{subequations}\label{KKT}
\begin{align}
&[(T_{\ell_k}^TT_{\ell_k})^{2} +\lambda I]w(\lambda)=-(T_{\ell_k}^TT_{\ell_k})^{1/2} \|g_k\|e_1,\label{KKTy}\\
&\lambda(\Delta_k-\|w(\lambda)\|)=0,\label{KKTy2}\\
&\lambda\geq 0,\\
&\|w(\lambda)\|\leq\Delta_k.
\end{align}
\end{subequations}

To solve \eqref{KKT} we look for a couple $(w(\lambda_k),\lambda_k)$ satisfying \eqref{KKTy} and such that $\|z(\lambda_k)\|=\|w(\lambda_k)\|=\Delta_k$ and we employ Newton's method to compute such $\lambda_k$, cf. \cite[ch. 7]{cgt}. We will see that, as in \cite{ellittica}, this is exactly what is required, as we will show that $\lambda_k>0$  and by \eqref{KKTy2} it follows $\|w(\lambda_k)\|=\Delta_k$. At each iteration of Newton's method we need to solve a sequence of linear systems of the form \eqref{KKTy} for fixed $\lambda$. The linear systems \eqref{KKTy} can be solved directly, due to the small size of the matrix. 

Finally,  given $(w_k(\lambda_k),\lambda_k)$ solution to \eqref{KKT}, $z(\lambda_k)$ can be recovered as $z(\lambda_k)=Q_{\ell_k} w(\lambda_k)$ and
 $p(\lambda)$ defined in \eqref{plambda}  can be approximated by \eqref{sqrt} with $b=z(\lambda)$, so that
\begin{subequations}
\begin{align}
p(\lambda)&=\ql(\TT)^{1/2}\ql^Tz(\lambda)\label{def_pl},\\
p_k=p(\lambda_k)&=\ql(\TT)^{1/2}\ql^Tz(\lambda_k)\label{def_p}.
\end{align}
\end{subequations} 
We remark again that, due to the small size of the matrix, $(\TT)^{1/2}$ can be computed directly, via an SVD decomposition.

	\section{Choice of the trust-region radius}\label{sec_radius}
	Taking into account that the problem is ill-posed, we aim at defining a trust-region method that shows regularizing properties. 
	To this end 
	we require that the step $p_k=p(\lambda_k)$ satisfies the following condition:
	\begin{equation}\label{qcond}
	\|\ql^Tm_k(p(\lambda_k))\|\geq q \|g_k\|,
	\end{equation}
	for $q$ a given constant in $(0,1)$.  We will show that this property  ensures regularizing properties to the method and that it can be satisfied by an appropriate choice of the trust-region radius.

	\begin{remark}
		This is a different condition from the so-called $q$-condition, employed in \cite{ellittica}:  
		\begin{equation}\label{qcond_ellittica}
		\|m_k(p(\lambda_k))\|\geq q\|g_k\|.
		\end{equation}
		 
		We will refer to \eqref{qcond} as the \textit{projected $q$-condition},  as $\|\ql\ql^T w\|=\|\ql^Tw\|$, for any $w\in \mathbb{R}^n$, and therefore condition \eqref{qcond} is equivalent to $\|\ql\ql^Tm_k(p(\lambda_k))\|\geq q \|g_k\|$.

	 Notice that if $\ell_k=n$, $\ql\ql^T$ is the identity matrix and the original condition \eqref{qcond_ellittica} is recovered.  
	\end{remark}
	We need to introduce this new condition to take into account that we are restricting the solution of \eqref{system_appr} to the subspace spanned by the columns of $Q_{\ell_k}$.
 	 
	We can prove that it exists a $\lambda_k^{q}$ such that  $\|\ql^Tm_k(p(\lambda_k^q))\|= q \|g_k\|$ 
	and that condition  \eqref{qcond} is satisfied for all $\lambda_k\geq\lambda_k^{q}$. 
	\begin{lemma}\label{lemma_lambdaq}
		$\|\ql^Tm_k(p(\lambda))\|$ is a monotone function of $\lambda$ such that 
		\begin{align*}
		&\lim_{\lambda\rightarrow\infty}\|\ql^Tm_k(p(\lambda))\|=\|g_k\|,\\
		&\lim_{\lambda\rightarrow 0}\|\ql^Tm_k(p(\lambda))\|=0.
		\end{align*}
		Moreover, it exists $\lambda_k^{q}>0$ such that 
		\begin{align}
		\lambda_k^{q}\leq \frac{q}{1-q}\|\TT\|^2,\label{bound}\\
		\|\ql^Tm_k(p(\lambda_k^{q}))\|=q\|g_k\|,\label{uguale}
		\end{align}
		and condition \eqref{qcond} is satisfied for all $\lambda_k\geq\lambda_k^{q}$.
	\end{lemma}
	\proof{
		We first prove that $\|\ql^Tm_k(p(\lambda))\|$ is a monotonic increasing function of $\lambda$. 
		Let $U\Sigma U^T$ be the $\mathrm{SVD}$  decomposition of matrix $\tl^T\tl$. Then $(\tl^T\tl)^2=U\Sigma^2 U^T$ and by \eqref{diag1} $Q_{\ell_k}^T B_k Q_{\ell_k}=\tl^T\tl=U\Sigma U^T$. 
		From \eqref{KKTy} and \eqref{def_p} it follows:
		\begin{align}
		w(\lambda)&=-U(\Sigma^2+\lambda I)^{-1}\Sigma^{1/2} U^T\ql^Tg_k,\label{w}\\
		z(\lambda)&=-\ql U(\Sigma^2+\lambda I)^{-1}\Sigma^{1/2} U^T\ql^Tg_k,\label{z}\\
		p(\lambda)&=\ql (\A)^{1/2} \ql^Tz(\lambda)= -\ql U\Sigma^{1/2}(\Sigma^2+\lambda I)^{-1}\Sigma^{1/2} U^T\ql^Tg_k\label{p}.
		\end{align}
		Moreover by \eqref{diag1}
		\begin{align*}
		\ql^Tm_k(p(\lambda))=&\ql^TB_kp(\lambda)+\ql^Tg_k=-\ql^TB_k\ql U\Sigma(\Sigma^2+\lambda I)^{-1} U^T\ql^Tg_k+\ql^Tg_k\\
		=&- U\Sigma^2(\Sigma^2+\lambda I)^{-1} U^T\ql^Tg_k+U(\Sigma^2+\lambda I)(\Sigma^2+\lambda I)^{-1}U^T\ql^Tg_k\\
		=&\lambda U (\Sigma^2+\lambda I)^{-1}U^T\ql^Tg_k.
		\end{align*}
		Then, denoting $r=U^T\ql^Tg_k$ and $\sigma_1,\dots,\varsigma^k_{\ell_k}$ the singular values of $\TT$, we obtain  
		\begin{equation}\label{norm_mod}
		\|\ql^Tm_k(p(\lambda))\|^2=\sum_{i=1}^{\ell_k}\left(\frac{\lambda }{\sigma_i^2+\lambda}\right)^2r_i^2
		\end{equation}
		and taking the derivative with respect to $\lambda$ 
		$$
		\frac{d}{d\lambda}\|\ql^Tm_k(p(\lambda))\|^2=\sum_{i=1}^{\ell_k}\frac{2\lambda}{\sigma_i^2+\lambda}\frac{\sigma_i^2 }{(\sigma_i^2+\lambda)^2}r_i^2>0.
		$$
		Moreover from \eqref{norm_mod}
		\begin{align*}
		\lim_{\lambda\rightarrow\infty}\|\ql^Tm_k(p(\lambda))\|&=\|r\|=\|g_k\|\\
		\lim_{\lambda\rightarrow 0}\|\ql^Tm_k(p(\lambda))\|&=0.
		\end{align*}
		Then, $\|\ql^Tm_k(p(\lambda))\|$ is a monotone function of $\lambda$ that varies between 0 and $\|g_k\|$.
		
		 Then, as $0<q<1$, there exists $\lambda_k^{q}$ such that \eqref{uguale} holds and $\lambda_k^{q}>0$. Due to the monotonicity of $\|\ql^Tm_k(p(\lambda_k^{q}))\|$, condition \eqref{qcond} is satisfied for all $\lambda_k\geq\lambda_k^{q}$.	

	Finally, we derive the bound \eqref{bound} on $\lambda_k^{q}$. It holds, by \eqref{norm_mod}, that
		\begin{eqnarray*}
		\|\ql^Tm_k(p(\lambda_k^{q}))\|^2&\geq& \left(\frac{\lambda_k^{q}}{\lambda_k^{q}+\|\TT\|^2}\right)^2\|g_k\|^2.
		\end{eqnarray*}
			Then, since $\lambda_k^{q}$ satisfies \eqref{uguale}, it holds
			$$
			q\|g_k\|=\|\ql^Tm_k(p(\lambda_k^{q}))\|\geq \frac{\lambda_k^{q}}{\lambda_k^{q}+\|\TT\|^2}\|g_k\|
			$$
			and  \eqref{bound} follows.
					
			}
	
	This lemma has the following important consequence. 
	
	\begin{lemma}\label{lemma_qcond}
		If $p(\lambda_k)$ given in \eqref{p} is such that the projected $q$-condition \eqref{qcond} is satisfied, then the trust-region constraint $\|w(\lambda_k)\|=\Delta_k$ is active.  
	\end{lemma}
	This is a consequence of the previous Lemma. Indeed,  as $\lambda_k\geq\lambda_k^q>0$, from \eqref{KKTy2} the trust-region must be active. 
	
	In the following lemma we prove that a suitable choice of the trust-region radius provides a parameter $\lambda$ for which condition \eqref{qcond} is satisfied. This is really important, as it allows for the sought automatic rule to satisfy the desired projected $q$-condition \eqref{qcond}. 
	\begin{lemma}
		If 
		\begin{equation}\label{delta}
		\Delta_k\leq (1-q)\frac{\|(\A)^{1/2}\ql^Tg_k\|}{\|\A\|^2},
		\end{equation}
		 then the projected $q$-condition \eqref{qcond} is satisfied. 
	\end{lemma}
	\proof{
		From \eqref{KKTy}
		$$
		\|z(\lambda_k^{q})\|=\|\ql^Tw(\lambda_k^{q})\|=\|w(\lambda_k^{q})\|\geq \frac{\|(\A)^{1/2}\ql^Tg_k\|}{\|(\TT)^2+\lambda_k^{q}I\|}.
		$$
		From \eqref{bound}
		$$
		\|(\A)^2+\lambda_k^{q} I\|\leq \frac{1}{1-q}\|\A\|^2.
		$$
		By construction $\|z(\lambda_k)\|\leq \Delta_k$. Then, 
		\begin{equation}
		\|z(\lambda_k)\|\leq \Delta_k \leq  (1-q)\frac{\|(\A)^{1/2}\ql^Tg_k\|}{\|\A\|^2}\leq\frac{\|(\A)^{1/2}\ql^Tg_k\|}{\|(\A)^2+\lambda_k^{q} I\|}\leq \|z(\lambda_k^{q})\|.
		\end{equation}
		From \eqref{z}, 
		$$\|z(\lambda)\|^2=\sum_{i=1}^{\ell_k}\frac{\sigma_i}{(\sigma_i^2+\lambda)^2}r_i^2$$ and therefore $\|z(\lambda)\|$
		 is a decreasing function of $\lambda$. 
		Then $\lambda_k\geq \lambda_k^{q}$ and from Lemma \ref{lemma_lambdaq} the projected $q$-condition is satisfied. 
	}

The resulting regularizing inexact Tikhonov method we propose, in its elliptical trust-region implementation, is sketched in Algorithm \ref{algoTR}. We underline that we initialize 
$q_1=g_k$ in GKLB and that 
the step acceptance criterion in Step 4.3 
 is different from that given in Step 2.3 of Algorithm \ref{algoSVD}.  Indeed, we compare the actual reduction in function values with that  predicted by the model
 \begin{equation}\label{model_phi} 
  \Phi_k(w)=\frac{1}{2}w^T (T_{\ell_k}^TT_{\ell_k})^2 w +w^T Q_{\ell_k}^T \tilde s_k+f_{\delta}(\xkd)
  \end{equation} 
 at $w_k:=w(\lambda_k)$,  rather than the classical model $\frac{1}{2}p^TB_kp+ p^Tg_k+f_{\delta}(\xkd)$. 
  This is necessary to take into account the sources of inexactness introduced in the method and to prove the theoretical results in the following section. 
\begin{algorithm}
	\caption{$k$-th step of the regularizing hybrid Lanczos trust-region method}
	\label{algoTR}
	Given $\xkd$,    $\eta\in (0,1)$, $\gamma\in (0,1)$, $0<C_{\min}<C_{\max}$, $q\in (0,1)$. 
	\vskip 1pt
	Exact data:  $y$; Noisy data:  $y^{\delta}$.
	\begin{description}
		\item[1.] Choose $1\leq \ell_k\leq n$.
		\item[2.] Compute $\tilde{s}_k$ as in \eqref{tsk}, with $Q_{\ell_k}$ and $T_{\ell_k}$ obtained from GKLB($\Jdk,g_k,\ell_k)$ (Algorithm \ref{algo_lanczos}).
		\item[3.]   Choose  $\displaystyle{\Delta_k}\in \left.\left [C_{\min}\|\tilde{s}_k\|, \, \min\left\{ C_{\max},  \frac{1-q}{\|\TT\|^2}  \right\}\|\tilde{s}_k\|\right]\right.$.
		\item[4.]   Repeat\\ 
		4.1 Compute $(T_{\ell_k}^TT_{\ell_k})^2$ and find the pair $(w_k,\lambda_k)$ solution of KKT conditions \eqref{KKT}. \\
		4.2 Set $z_k=Q_{\ell_k} w_k$ and $p_k=Q_{\ell_k} (\TT)^{1/2}Q_{\ell_k}^Tz_k$.\\
		4.3  Set    $\Phi_k(w_k)=\frac{1}{2}w_k^T (T_{\ell_k}^TT_{\ell_k})^2 w_k +w_k^T Q_{\ell_k} \tilde s_k+f_{\delta}(\xkd)$ and compute  
		$${\pi_k(w_k)}=\frac{f_\delta(\xkd)-f_\delta(\xkd+p_k)}{f_\delta(\xkd)-\Phi_k(w_k)}$$ 
		4.4 If $\pi_k(w_k)< \eta$, set  ${\Delta_{k}=\gamma \Delta_k}$.\\
		\hspace*{-20pt}		Until $\pi_k(w_k)\ge \eta$.
		\item [5.] Set  ${x_{k+1}^\delta=\xkd+p_k}$.
	\end{description}
\end{algorithm}

	We prove a technical lemma, which will be useful in the following. 

\begin{lemma}\label{lemma_T}
Let $e_1$ be first vector of the canonical basis. Let $J(\xkd)$ be of rank $r$ and let $\varsigma^k_1\geq\dots\geq\varsigma^k_{r}>0$ be the nonzero singular values of $J(\xkd)$. Then, 
at each iteration of Algorithm \ref{algoTR} the following inequality holds:
$$\tilde{t}_k:=e_1^T\tl^T\tl e_1\geq(\varsigma^k_{r})^2.$$\end{lemma}
\proof{
From Algorithm \ref{algo_lanczos} it holds $e_1^T \tl^T\tl e_1=\alpha_1^2=\frac{\|J(\xkd)g_k\|^2}{\|g_k\|^2}=\frac{\|J(\xkd)J(\xkd)^TF(\xkd)\|^2}{\|J(\xkd)^TF(\xkd)\|^2}$.
Let $J(\xkd)=\bar{U}\bar{\Sigma}\bar{ V}^T$ be the singular value decomposition of $J(\xkd)$ and $r=\bar{U}^TF(\xkd)$.
The thesis follows from
\begin{equation*}
\|J(\xkd)J(\xkd)^TF(\xkd)\|^2=\|\bar{U}\bar{\Sigma}\bar{\Sigma}^T\bar{U}^TF(\xkd)\|^2=\|\bar\Sigma\bar\Sigma^Tr\|^2=\sum_{i=1}^{r}\varsigma_i^4 r_i^2\geq \varsigma_r^2\sum_{i=1}^{r}\varsigma_i^2 r_i^2,
\end{equation*}
and
\begin{equation*}
\|J(\xkd)^TF(\xkd)\| ^2=\| \bar V\bar \Sigma^T\bar U^TF(\xkd)\|^2=\|\bar\Sigma^T r\|^2=\sum_{i=1}^{r}\varsigma_i^2 r_i^2.
\end{equation*}
}

	\section{Theoretical analysis: convergence and complexity}\label{sec_algo}
In this section we will discuss the convergence, the complexity and the computational cost of the method in Algorithm \ref{algoTR}. To prove these results, we make the following assumptions on $J$.
 \begin{assumption}\label{hp_bound}
	Assume that for all $x$ in a neighbourhood of the level set $\mathcal{L}=\{x\in\mathbb{R}^n\;\; s.t. \;\; f_\delta(x)\leq f_\delta(x_0)\}$,  there exists $K_J>0$ such that $\|J(x)\|\leq K_J$.
\end{assumption}

\begin{assumption}\label{hp_lip}
	The gradient of $f$ is Lipschitz continuous in a neighbourhood of the level set $\mathcal{L}=\{x\in\mathbb{R}^n\;\; s.t. \;\; f_\delta(x)\leq f_\delta(x_0)\}$   with Lipschitz constant $L$.
\end{assumption}

\begin{assumption}\label{hp_rank}
Assume that for all $x$ in a neighbourhood of the level set $\mathcal{L}=\{x\in\mathbb{R}^n\;\; s.t. \;\; f_\delta(x)\leq f_\delta(x_0)\}$
$J(x)$ has full rank.

\end{assumption}

 We remark that this last assumption makes the proof of the convergence properties easier, but it is not necessary. We will discuss the properties and implementation of the method in case of rank deficient Jacobian in Section \ref{sec_analysis}.

 We observe that  the choice  $q_1=g_k$ in GKLB Algorithm by \eqref{tsk} yields
\begin{equation} \label{phi0}
\nabla \Phi_k(0)=\ql^T\tilde s_k=  (\tl^T\tl)^{1/2} \ql^Tg_k=\|g_k\| (\tl^T\tl)^{1/2} e_1.
\end{equation}
This is a crucial property and will be used in the following lemma to quantify the decrease provided by the so-called Cauchy step and to prove global convergence of the method  in the noise-free case  and complexity results in presence of noise.
 The Cauchy step $w_k^c$ is the minimizer of $\Phi_k(w)$ along $-\nabla \Phi_k(0)$ within the trust region, i.e. it is the vector $w_k^c:=-\alpha_k\nabla \Phi_k(0) $, where $\alpha_k$ satisfies
\begin{equation}\label{alpha_cauchy}
\alpha_k=\arg\min_{\alpha: \,\|\alpha \nabla \Phi_k(0)\|\leq \Delta_k}\Phi_k(-\alpha \nabla \Phi_k(0)).
\end{equation}
\begin{lemma}\label{lemma_cauchy}
The decrease in the model $\Phi_k(w)$ given in \eqref{model_phi}, achieved by the Cauchy step $w_k^c$, is such that
\begin{equation}\label{cd}
\Phi_k(0)-\Phi_k(w_k^c)\geq
\frac{1}{2}\tilde{t}_k^{1/2} \|g_k\| \min \Bigg\{\Delta_k,\frac{\tilde{t}_k^{1/2} \|g_k\|}{ \|T_{\ell_k}^TT_{\ell_k}\|^2}\Bigg\},
\end{equation}
with $\tilde{t}_k$ defined in Lemma \ref{lemma_T}.
\end{lemma}
\proof{
From Lemma 4.3 of \cite{nw} it follows
$$
\Phi_k(0)-\Phi_k(w_k^c)\geq
\frac{1}{2}\|\nabla \Phi_k(0)\| \min \Bigg\{\Delta_k,\frac{ \|\nabla \Phi_k(0)\|}{ \|T_{\ell_k}^TT_{\ell_k}\|^2 }\Bigg\}.
$$
Then the thesis follows 
taking into account that by \eqref{phi0} we have 
$$
\|\nabla \Phi_k(0)\|=\tilde t_k^{1/2} \|g_k\|.
$$
}

Let us now present a lower bound on the trust-region radius at $k$-iteration. This will be used to prove the complexity result and also ensures that the \textit{repeat} cycle at Step 4  of Algorithm \ref{algoTR} terminates in a finite number of steps.

\begin{lemma}\label{Lemma_bound_delta}
Assume Assumptions \ref{hp_bound}-\ref{hp_rank} hold. 
In Algorithm \ref{algoTR} the trust region radius $\Delta_k\geq\Delta_{\min,k}$, with $\Delta_{\min,k}:=\gamma(1-\eta)\frac{\varsigma^k_{\ell_k}\|g_k\|}{K_J^2\left(L+ K_J^2\right)}$.
\end{lemma}
\proof{
From Step 4.3 of Algorithm \ref{algoTR} 
\begin{equation*}
\pi_k(w_k)-1=-\frac{f_\delta(\xkd+p_k)-\Phi_k(w_k)}{f_\delta(\xkd)-\Phi_k(w_k)}.
\end{equation*}
We remark that from \eqref{diag1} it holds $\|\tl^T\tl\|\leq \|B_k\|$. From this, 
Assumptions \ref{hp_bound} and \ref{hp_lip}, the definition   of   $z_k=\ql w_k$, the definition of $\Phi_k(w)$ given in \eqref{model_phi} and relations \eqref{tsk} and  \eqref{def_p} it holds (reminding also \eqref{notation}):
\begin{align*}
\lvert f_\delta(\xkd+p_k)-\Phi_k(w_k)\rvert &=\Bigg\lvert \int_{0}^{1}\nabla f_\delta(\xkd+tp_k)^Tp_k\,dt-\frac{1}{2}w_k^T(T_{\ell_k}^TT_{\ell_k})^2w_k-z_k^T\tilde s_k\Bigg\rvert\\
&\stackrel{\eqref{tsk}}{=}\Bigg\lvert \int_{0}^{1}\nabla f_\delta(\xkd+tp_k)^Tp_k\,dt-\frac{1}{2}w_k^T(T_{\ell_k}^TT_{\ell_k})^2w_k-z_k^TQ_{\ell_k} (\tl^T\tl)^{1/2}Q_{\ell_k}^Tg_k\Bigg\rvert\\
&\stackrel{\eqref{def_p}+\eqref{notation}}{=}\Bigg\lvert \int_{0}^{1}[\nabla f_\delta(\xkd+tp_k)-\nabla f_\delta(\xkd)]^Tp_k\,dt-\frac{1}{2}w_k^T(T_{\ell_k}^TT_{\ell_k})^2w_k\Bigg\rvert\\
&\leq \frac{L}{2}\|p_k\|^2+ \frac{1}{2}K_J^4\|w_k\|^2\leq \left(\frac{L}{2}\|\tl^T\tl\|+ \frac{1}{2}K_J^4\right)\|z_k\|^2\\
&\leq \left(\frac{L}{2}K_J^2+ \frac{1}{2}K_J^4\right)\Delta_k^2=\frac{K_J^2}{2}(L+ K_J^2)\Delta_k^2.
\end{align*}
Moreover, since $J(\xkd) $ has full rank by hypothesis, Lemma \ref{lemma_T}  guarantees that $\tilde t_k>0$.
Assume that   $\Delta_k \le \frac{\tilde{t}_k^{1/2} \|g_k\|}{ \|T_{\ell_k}^TT_{\ell_k}\|^2}$. 
Then, since $w_k$ is the solution of \eqref{min_proj},
from Lemma \ref{lemma_T} and \ref{lemma_cauchy} it follows
\begin{equation*}
f_\delta(\xkd)-\Phi_k(w_k)\geq f_\delta(\xkd)-\Phi_k(w_k^c)\geq \frac{1}{2}\varsigma^k_{\ell_k}\|g_k\| \Delta_k.
\end{equation*}
Then, 
\begin{equation*}
\lvert\pi_k(w_k)-1\rvert=\Bigg\lvert\frac{f_\delta(\xkd+p_k)-\Phi_k(w_k)}{f_\delta(\xkd)-\Phi_k(w_k)}\Bigg\rvert\leq \frac{K_J^2\left(L+ K_J^2\right)\Delta_k}{\|g_k\| \varsigma^k_{\ell_k}}.
\end{equation*}
Then, $\pi_k(w_k)\geq\eta$ and the iteration is successful whenever 
 $$
 \Delta_k \le\frac{(1-\eta)\varsigma^k_{\ell_k} \|g_k\| }{K_J^2\left(L+ K_J^2\right)}=
  \min\Bigg\{ \frac{\varsigma^k_{\ell_k} \|g_k\|}{ \|T_{\ell_k}^TT_{\ell_k}\|^2}, (1-\eta)\frac{ \varsigma^k_{\ell_k}\|g_k\|}{K_J^2\left(L+ K_J^2\right)}\Bigg\}.
 $$
From the updating rule of $\Delta_k$ at Step 4 of Algorithm \ref{algoTR} it must hold $\Delta_k\geq \Delta_{\min,k}$.
}

We are now ready to state the complexity result.

\begin{theorem}\label{teo_complexity}
Let Assumptions \ref{hp_bound}-\ref{hp_rank} hold. Assume that  there exists   
$\varsigma>0$ such that $\varsigma^k_{\ell_k}>\varsigma$ for any $k$ and  $f_\delta$ is bounded below on the level set $\mathcal{L}=\{x\in\mathbb{R}^n \text{ s.t. } f(x)\leq f(x_0)\}$.
Given a positive constant $\epsilon>0$, the method in Algorithm \ref{algoTR} takes at most $O(\epsilon^{-2})$ iterations to achieve $\|g_k\|\leq \epsilon$.
\end{theorem}
\proof{
Let $k_\epsilon$ be the first iteration index such that $\|g_k\|<\epsilon$ and $f_{\delta}^{\min}$ be the lower bound of $f_\delta$ in $\mathcal{L}$.
 From Step 4 of Algorithm \ref{algoTR}, Lemma \ref{lemma_T}, Lemma \ref{Lemma_bound_delta}  and \eqref{cd} it holds:
\begin{align*}
f_\delta(x_0)-f_{\delta}^{\min}&\geq f_\delta(x_0)-f_\delta(x_{k_\epsilon}^\delta)=\sum_{k=0}^{\ke} f_\delta(\xkd)-f_\delta(\xkd+p_k)\geq \eta \sum_{k=0}^{\ke} f_\delta(\xkd)-\Phi_k(w_k)\\
&\geq \frac{\ke \eta}{2}\epsilon \varsigma^k_{\ell_k}\min \Bigg\{\Delta_{\min,k}, \frac{ \|\nabla \Phi_k(0)\|}{ \|T_{\ell_k}^TT_{\ell_k}\|^2 }\Bigg\}\\
&\geq \frac{\ke \eta}{2}\epsilon  \varsigma\gamma(1-\eta)\frac{\varsigma\epsilon }{K_J^2\left(L+ K_J^2\right)}.
\end{align*}
 Then, $\ke\leq \left \lceil{ \frac{f_{\delta}(x_0)-f_{\delta}^{\min}}{C\epsilon^2}}\right \rceil $ where $C=\left(\frac{ \gamma\eta(1-\eta) \varsigma^2 }{ 2 K_J^2(L+ K_J^2)}\right)^{-1}$.
 
}

 In the case of noisy data, the process is stopped before the convergence is reached, to avoid approaching a solution of the noisy problem, according to \eqref{discrepance}. 
Theorem \ref{teo_complexity} guarantees that the stopping criterion is satisfied 
after $O(\delta^{-2})$ iterations.  
 From this theorem, the following corollary easily follows, which states the global convergence of the proposed method. 
\begin{corollary}
	Let Assumptions in Theorem \ref{teo_complexity}  hold. 	The sequence $\{x_k\}$ generated by  Algorithm \ref{algoTR} satisfies
	\begin{equation*}
	\liminf_{k\rightarrow \infty} \nabla f(x_k)=0.
	\end{equation*} 
\end{corollary}

In terms of computational cost, assuming $\ell_k=\ell$, for all $k$'s, in the worst case the procedure requires $O(\delta^{-2}\ell)$ matrix-vector products.

\section{Theoretical analysis: regularizing properties}\label{sec_analysis}
In this section we prove that it is possible to control how the distance of the current solution approximation from an exact solution of the problem changes from an iteration to the other, and that this decreases if the dimension of the Lanczos space is large enough. We first state all the assumptions we make. They are the same as in \cite{ellittica}. 

\begin{assumption}\label{ass_taylor}
	Given $ x,\tilde x$ in a suitable neighbourhood of the solution $x^\dagger$ to \eqref{problem}, the following inequality holds:
	\begin{equation}\label{taylor1}
	\| \nabla f(\tilde{x})-\nabla f(x)-J(x)^TJ(x)(\tilde{x}-x)\|\leq
	(c \| \tilde{x}-x\|+\sigma) \| \nabla f(x)-\nabla f(\tilde{x})\|,
	\end{equation}
	where $\nabla f$ is the gradient of $f$ and $J$ is the Jacobian of matrix of $F$, for suitable constants  $c>0$ and $\sigma\in (0,q)$.  
\end{assumption} 

As explained in \cite{ellittica}, the above assumption is reasonable if $\sigma$ is interpreted as a bound for $\|S(x)\|$, with $S(x)=\sum_{i=1}^{m}(F_i(x)-y)\nabla^2 F_i(x)$, the term containing the  second-order information in $\nabla^2f(x)$. The constant $\sigma$ can therefore be interpreted as 
 a combined measure of the nonlinearity and residual size of the problem.

We then prove  this technical lemma.
\begin{lemma}
Let $m_k(p(\lambda))$ be defined in \eqref{def_mod}. Let $Q_{\ell_k}, T_{\ell_k}$ be the output of GKLB($\Jdk$,$g_k$,$\ell_k$). Let $w(\lambda)$ be such that \eqref{KKTy} hold for $\lambda>0$ and let $z(\lambda)=Q_{\ell_k}w(\lambda)$ and $p(\lambda)$ given by \eqref{def_pl}.
Then,
\begin{equation}\label{rel_mod_p}
-\lambda Q_{\ell_k}^Tp(\lambda)=\TT\ql^T m_k(p(\lambda)).
\end{equation}
\end{lemma}
\proof{
From \eqref{def_pl}, \eqref{KKTy} and  the definition of $w(\lambda)$  it follows that 
$$
Q_{\ell_k}^Tp(\lambda)=(\TT)^{1/2}w(\lambda)= -(\TT)^{1/2}[(T_{\ell_k}^TT_{\ell_k})^{2} +\lambda I]^{-1}(T_{\ell_k}^TT_{\ell_k})^{1/2} Q_{\ell_k}^Tg_k.
$$
Then
$$
(T_{\ell_k}^TT_{\ell_k})^{-1/2}[(T_{\ell_k}^TT_{\ell_k})^{2} +\lambda I](\TT)^{-1/2}Q_{\ell_k}^Tp(\lambda)=- Q_{\ell_k}^Tg_k
$$
and  
$$
T_{\ell_k}^TT_{\ell_k}Q_{\ell_k}^Tp(\lambda)+Q_{\ell_k}^Tg_k=-\lambda(\TT)^{-1}Q_{\ell_k}^Tp(\lambda).
$$
Multiplying the previous relation by $\TT$, we obtain
\begin{equation}\label{rel1}
(T_{\ell_k}^TT_{\ell_k})^2Q_{\ell_k}^Tp(\lambda)+\TT Q_{\ell_k}^Tg_k=-\lambda Q_{\ell_k}^Tp(\lambda).
\end{equation}
Moreover, multiplying both sides of \eqref{def_mod} by $\TT Q_{\ell_k}^T$, 
$$
\TT Q_{\ell_k}^Tm_k(p(\lambda))=\TT Q_{\ell_k}^TB_kp(\lambda)+\TT Q_{\ell_k}^Tg_k.
$$
Finally, 
\begin{align*}
\TT Q_{\ell_k}^Tm_k(p(\lambda))&\stackrel{\eqref{def_pl}}{=}\TT Q_{\ell_k}^TB_kQ_{\ell_k}(\TT)^{1/2}Q_{\ell_k}^Tz(\lambda)+\TT Q_{\ell_k}^Tg_k\\
&\stackrel{\eqref{diag1}}{=}(\TT)^2Q_{\ell_k}^TQ_{\ell_k}(\TT)^{1/2}Q_{\ell_k}^Tz(\lambda)+ \TT Q_{\ell_k}^Tg_k\\
&\stackrel{\eqref{def_p}}{=}(\TT)^2Q_{\ell_k}^Tp(\lambda)+\TT Q_{\ell_k}^Tg_k\stackrel{\eqref{rel1}}{=}-\lambda Q_{\ell_k}^Tp(\lambda).
\end{align*}
}

\begin{remark}\label{remark_rank1}
 We highlight that the structure of the step is fundamental to prove \eqref{rel_mod_p}.  If we compute $w_k$ from \eqref{sub_w} rather than from  \eqref{min_proj}, i.e. if we keep the rank-one perturbation in \eqref{rank1}, then the crucial relation \eqref{rel_mod_p}  between the step and the model does not hold. Note that, this condition is guaranteed to be satisfied in the exact case, i.e.  when   $\ell_k=n$ is taken.  
\end{remark}
In the following lemma we provide a condition ensuring that  the distance between the computed solution and the true solution  decreases from iteration $k$ to iteration $k+1$, provided that $x_k$ is sufficiently close to the true solution and $\|Q_{\ell_k} T_{\ell_k}^TT_{\ell_k} Q_{\ell_k}^T-\ql\ql^TB_k\|$ is sufficiently small.

\begin{lemma}\label{monotone_decay}
	Assume that $x^\dagger$ is a solution of \eqref{problem}. Let $e_k=x^\dagger-\xkd$, $x_{k+1}^\delta=\xkd+p_k$ with $p_k=p(\lambda_k)$ defined in \eqref{def_p}, and 	$m_k(p(\lambda))$ defined  in \eqref{def_mod}. Assume that 
	 there exists $\theta_k>1$ such that the following condition holds:
	  \begin{equation}\label{condition_theta}
	 \| \ql^Tm_k(e_k)\|\leq \frac{1}{\theta_k}\|\ql^Tm_k(p_k)\|,\;\;\; \theta_k>1.
	 \end{equation}  
	 Then
	\begin{align}
	\|x_{k+1}^\delta-x^\dagger\|^2-\|\xkd-x^\dagger\|^2\leq&\frac{2}{\lambda_k} \left(\frac{1}{\theta_k}-1\right) \| \ql^Tm_k(p_k)\|^2
	 \nonumber\\
	+&\| \ql^Tm_k(p_k)\|\|Q_{\ell_k} T_{\ell_k}^TT_{\ell_k} Q_{\ell_k}^T-\ql\ql^TB_k\|\|x_{k}^\delta-x^\dagger\|.
	\label{decr_errore}
	\end{align}
\end{lemma}

\proof{
	Note that
	\begin{align}
	\|x_{k+1}^\delta-x^\dagger\|^2 -\|\xkd-x^\dagger\|^2  =&
	2\langle x_{k+1}^\delta-\xkd,\xkd-x^\dagger\rangle +\|x_{k+1}^\delta-\xkd\|^2\nonumber\\
	=&-2\langle p_k,e_k\rangle +\|p_k\|^2\label{inizio}.
	\end{align}

Remark that from \eqref{def_p}, $p_k\in\mathcal{R}(Q_{\ell_k})$, so that $Q_{\ell_k} Q_{\ell_k}^T p_k=p_k$. From this and \eqref{rel_mod_p} it holds:
\begin{align}
\langle p_k,e_k\rangle&=\langle Q_{\ell_k} Q_{\ell_k}^Tp_k,e_k\rangle=\langle Q_{\ell_k}^Tp_k,Q_{\ell_k}^Te_k\rangle=-\frac{1}{\lambda_k}\langle  Q_{\ell_k}^Tm_k(p_k),T_{\ell_k}^TT_{\ell_k}\ql^T e_k\rangle\nonumber\\
&=-\frac{1}{\lambda_k}\langle Q_{\ell_k}^Tm_k(p_k),( T_{\ell_k}^TT_{\ell_k} Q_{\ell_k}^T-Q_{\ell_k}^TB_k)e_k\rangle-\frac{1}{\lambda_k}\langle Q_{\ell_k}^Tm_k(p_k),Q_{\ell_k}^TB_ke_k\rangle\label{rel0}.
\end{align}

Moreover,
\begin{align}
\langle Q_{\ell_k}^Tm_k(p_k), Q_{\ell_k}^TB_k e_k\rangle= & \langle Q_{\ell_k}^Tm_k(p_k),Q_{\ell_k}^T(B_k e_k+g_k)\rangle -\langle  Q_{\ell_k}^Tm_k( p_k),Q_{\ell_k}^T(B_k p_k+g_k)\rangle +\langle Q_{\ell_k}^Tm_k(p_k),Q_{\ell_k}^TB_k p_k\rangle\nonumber\\
= & \langle Q_{\ell_k}^Tm_k(p_k),Q_{\ell_k}^Tm_k(e_k)\rangle -\langle Q_{\ell_k}^Tm_k(p_k),Q_{\ell_k}^Tm_k(p_k)\rangle +\langle Q_{\ell_k}^Tm_k(p_k),Q_{\ell_k}^TB_kp_k\rangle\label{rel}.
\end{align}
Again from \eqref{rel_mod_p}, \eqref{diag1} and the fact that $p_k\in\mathcal{R}(Q_{\ell_k})$,
\begin{align}
\langle p_k,p_k\rangle=&\langle \ql\ql^Tp_k,p_k\rangle=\langle \ql^Tp_k,\ql^Tp_k\rangle=-\frac{1}{\lambda_k}\langle \tl^T\tl\ql^Tm_k(p_k),\ql^Tp_k\rangle\nonumber\\
=&-\frac{1}{\lambda_k}\langle \ql^Tm_k(p_k),\tl^T\tl\ql^Tp_k\rangle=-\frac{1}{\lambda_k}\langle \ql^Tm_k(p_k),\ql^TB_kp_k\rangle\label{norm_p}.
\end{align}
Putting together \eqref{rel0}, \eqref{rel} and \eqref{norm_p}, we obtain:
\begin{align*}
\langle p_k,e_k\rangle=&-\frac{1}{\lambda_k}\langle \ql^Tm_k(p_k),(T_{\ell_k}^TT_{\ell_k} Q_{\ell_k}^T-\ql^TB_k)e_k\rangle-\frac{1}{\lambda_k}\langle \ql^Tm_k(p_k),\ql^Tm_k(e_k)\rangle \\
&+\frac{1}{\lambda_k}\langle \ql^Tm_k(p_k),\ql^Tm_k(p_k)\rangle +\langle p_k,p_k\rangle.
\end{align*}

From \eqref{inizio}, using \eqref{condition_theta}, it follows
\begin{align*}
\|x_{k+1}^\delta-x^\dagger\|^2-\|\xkd-x^\dagger\|^2
& \leq \frac{2}{\lambda_k}\|\ql^T m_k(p_k)\| \| \ql^Tm_k(e_k)\|-\frac{2}{\lambda_k}\|\ql^T m_k(p_k)\|^2-\|p_k\|^2\\
&+\frac{2}{\lambda_k}\|\ql^T m_k(p_k)\|\|Q_{\ell_k} T_{\ell_k}^TT_{\ell_k} Q_{\ell_k}^T-\ql\ql^TB_k\|\|e_k\|\\
&\leq\frac{2}{\lambda_k} \frac{1}{\theta_k}  \| \ql^Tm_k(p_k)\|^2-\frac{2}{\lambda_k}\| \ql^Tm_k(p_k)\|^2 -\|p_k\|^2\\
&+\frac{2}{\lambda_k}\|\ql^T m_k(p_k)\|\|Q_{\ell_k} T_{\ell_k}^TT_{\ell_k} Q_{\ell_k}^T-\ql\ql^TB_k\|\|e_k\|\\
&\leq\frac{2}{\lambda_k} \left(\frac{1}{\theta_k}-1\right) \| \ql^Tm_k(p_k)\|^2+ \\
&+\frac{2}{\lambda_k}\| \ql^Tm_k(p_k)\|\|Q_{\ell_k} T_{\ell_k}^TT_{\ell_k} Q_{\ell_k}^T-\ql\ql^TB_k\|\|e_k\|,
\end{align*}
which yields the thesis. 

}

 We now prove that assumption \eqref{condition_theta} in the previous lemma can be satisfied if the current approximation is close enough to $x^\dagger$. We first focus on the noise-free case.

	\begin{lemma} \label{lemma_q}
		Assume that $x^\dagger$ is a solution of \eqref{problem} and that $\delta=0$. 
		Assume further that there exist $\rho>0$, $c>0$  and $\sigma \in (0,q)$ such that Assumption \ref{ass_taylor} holds for any $ x,\tilde x \in \mathcal{B}_{2\rho} (x_k)$, and that  		\begin{equation}
			\|x_k-x^\dagger\|< \min\left\{\frac{q-\sigma}{ c}, \rho\right \} \label{locTR1}.
			\end{equation}
		Then,	it exists $\theta_k>1$ such that 
		\eqref{condition_theta} is satisfied. 
			\end{lemma}
	\proof{
		Applying  Assumption \ref{ass_taylor} with $\tilde{x}=x^\dagger$ and $x=x_{{k}}$ it holds $\nabla f(x^\dagger)=0$ and
		$$
		\|m_{{k}}(e_{{k}})\|=\| g_{{k}}+B_{{k}}(x^\dagger-x_{{k}})\|\leq (c \|e_{{k}}\|+\sigma) \| g_{{k}}\|.
		$$
		
		Then,  it holds 
		$$
		\|Q_{\ell_k}^Tm_{{k}}(e_{k})\|\leq \|m_{{k}}(e_{k})\|\leq (c \|e_{{k}}\|+\sigma)\|  g_{{k}}\|.
		$$
		If we let 
		$$
		\theta_{{k}}=\frac{q}{(c \|e_{{k}}\|+\sigma)}
		$$
		from the assumption it holds $\theta_{{k}}>1$  and 
		from the projected $q$-condition \eqref{qcond} we obtain
		$$
		\|Q_{\ell_k}^Tm_{{k}}(e_{{k}})\|\leq \frac{(c \|e_{{k}}\|+\sigma)}{q} \|Q_{\ell_k}^Tm_{{k}}(p_k)\|
		= \frac{1}{\theta_{{k}}} \|Q_{\ell_k}^Tm_{{k}}(p_{{k}})\| .
		$$
	}

We focus now on the noisy case.

\begin{lemma} \label{lemma_q_noise}
	Assume that $x^\dagger$ is a solution of \eqref{problem} and that $\delta>0$. Assume further that there exist $\rho>0$, $c>0$  and $\sigma \in (0,q)$ such that Assumption \ref{ass_taylor} holds for any $ x,\tilde x \in \mathcal{B}_{2\rho} (x_k)$, and that  	
	\begin{equation}
	\|\xkd-x^\dagger\|< \min\left\{\frac{\bar\tau(q-\sigma)-(1+\sigma)}{ (1+\bar\tau)c}, \rho\right \} \label{locTR2}, \quad \bar\tau>\frac{1+\sigma}{q-\sigma},
	\end{equation}
	with $\bar\tau$ defined in \eqref{discrepance}. 
	Then,	it exists $\theta_k>1$ such that 
	\eqref{condition_theta} is satisfied. 
\end{lemma}
\proof{
	Let us remind that $g_k:=\nabla f_\delta(x_k^\delta)$. 
	Applying  Assumption \ref{ass_taylor} with $\tilde{x}=x^\dagger$ and $x=x_{{k}}^\delta$ it holds $\nabla f(x^\dagger)=0$ and
	\begin{align*}
	\|m_{{k}}(e_{{k}})\|&=\| g_k+B_k(x^\dagger-\xkd)\|\leq \| \nabla f(\xkd)+J(\xkd)^T(y-y^\delta)+B_{{k}}(x^\dagger-\xkd)\|\\
	&\leq \|J(\xkd)\|\delta+(c \|e_{{k}}\|+\sigma) \| \nabla f(\xkd)\|\leq  \|J(\xkd)\|\delta+(c \|e_{{k}}\|+\sigma) \|g_k-J(\xkd)^T(y-y^\delta)\| \\
	&\leq (1+c \|e_{{k}}\|+\sigma) \|J(\xkd)\|\delta+(c \|e_{{k}}\|+\sigma) \| g_{{k}}\|.
	\end{align*}
	If the discrepancy principle \eqref{discrepance} is not satisfied, it follows that
\begin{align*}
\|m_{{k}}(e_{{k}})\|&\leq \left(\frac{(1+c \|e_{{k}}\|+\sigma) \|J(\xkd)\|}{\tau_k}+(c \|e_{{k}}\|+\sigma) \right)\| g_{{k}}\|=\left(\frac{(1+c \|e_{{k}}\|+\sigma) }{q\bar{\tau}}+\frac{c \|e_{{k}}\|+\sigma}{q} \right)q\| g_{{k}}\|.
\end{align*}
	If we let 
	$$
	\theta_{{k}}=\frac{\bar\tau q}{1+(1+\bar\tau)(c\|e_k\|+\sigma)}
	$$
	from the assumption it holds $\theta_{{k}}>1$  and 
	from the projected $q$-condition \eqref{qcond} we obtain
	$$
	\|Q_{\ell_k}^Tm_{{k}}(e_{{k}})\|\leq \|m_{{k}}(e_{{k}})\|\leq \frac{1}{\theta_{{k}}} \|Q_{\ell_k}^Tm_{{k}}(p_{{k}})\| .
	$$
}

\begin{remark}\label{remark_monotono}
	 The results in Lemmas \ref{monotone_decay}, \ref{lemma_q} and \ref{lemma_q_noise} allow to 
establish a relation between the error at two successive iterations. This allows then to control the distance of the current iterate from the sought solution. In particular, if $\|Q_{\ell_k} T_{\ell_k}^TT_{\ell_k} Q_{\ell_k}^T-\ql\ql^TB_k\|$ is small enough, the  norm of  the error decreases from an iteration to another.  	 	
	In fact, assuming $\|e_{k}\|<1$, we can rewrite \eqref{decr_errore} as
	\begin{eqnarray*}
	\|x_{k+1}^\delta-x^\dagger\|^2-\|x_{k}^\delta-x^\dagger\|^2\leq \|\ql^Tm_k(p_k)\|\left[\frac{2}{\lambda_k} \left(\frac{1}{\theta_k}-1\right)\| \ql^Tm_k(p_k)\|+\|Q_{\ell_k} T_{\ell_k}^TT_{\ell_k} Q_{\ell_k}^T-\ql\ql^TB_k\|\right].
	\end{eqnarray*}
	Then, recalling that $\theta_k>1$,  if \eqref{condition_theta} is satisfied
	and
	$
	\|Q_{\ell_k} T_{\ell_k}^TT_{\ell_k} Q_{\ell_k}^T-\ql\ql^TB_k\|
	$
	is sufficiently small, it follows $\|x_{k+1}^\delta-x^\dagger\|^2<\|x_{k}^\delta-x^\dagger\|^2$.
	We remark that the quantity $\|Q_{\ell_k} T_{\ell_k}^TT_{\ell_k} Q_{\ell_k}^T-\ql\ql^TB_k\|$ could be used as a measure for the choice of parameters $\ell_k$ in Algorithm \ref{algoTR}.  For $\ell_k=n$ this term is zero, so if $\ell_k=n$ for any $k$'s, we recover the theoretical results in \cite{ellittica}. 
\end{remark}
\begin{remark}\label{remark_singular}
  The assumption that  the Jacobian matrix has full rank in a neighborhood of the solution can be relaxed.  The method is well defined even if the Jacobian is rank deficient, provided that the  matrix $T_{\ell_k}^TT_{\ell_k}$ is invertible for all $k$'s. Indeed, the results in Section \ref{sec_analysis} 
still hold in the rank deficient case as only the inverse of  $T_{\ell_k}^TT_{\ell_k}$ is needed.  Moreover,  all the results in Section \ref{sec_algo} still hold if the Jacobian does not change rank asymptotically, as under this assumptions   $\varsigma^k_{r}>\varsigma$ for any $k$ and $\tilde t_k$ in Lemma 4 is bounded away from zero. 
We outline that in Algorithm GKLB we choose $q_1=g_k$. Then, $q_1$ belongs to the range space of $J(\xkd)^T$ and therefore Algorithm GKLB generates all vectors $q_\ell$  in the range space of $J(\xkd)^T$ (see Step 2.4). Consequently, if 
$J(\xkd)$ has rank $r$,  after at most $r$ iterations  the algorithm meets a null $\beta_{k+1}$. If the method stops at iteration $\ell<r$ it is possible to continue the process until $\ell=r$ 
choosing a new vector $p_{k+1}$ ($q_{k+1}$) which is orthogonal to the
previous $p_k$'s ($q_k$'s) \cite{golub1965}. 
  \end{remark}

\section{Numerical results}\label{sec_num}
 In this section we present the results of our numerical experimentation, aimed at showing that the method presented in this work in practice shares the same regularizing behaviour of the method in \cite{ellittica}, and that  employing   $\ell_k< n $ it provides a solution approximation in shorter CPU time without compromising the accuracy of the approximation found.
 
Two nonlinear ill-posed least-squares problems have been selected. 
{\bf Problem 6.1} arises from the discretization of a parameter identification problem, while {\bf Problem 6.2} is an image registration problem. 
In the following,  the Euclidean norm will be denoted by $\|\cdot\|$.
\begin{itemize}
\item {\bf Problem 6.1}: {\tt A 2D parameter identification problem.} The problem consists of reconstructing $c$ in the 2D-elliptic problem 
	\begin{subequations}\label{eq}
		\begin{gather}
		-\Delta u +c u  = \varphi  \text{ in } \Omega\label{eq1}\\
		u  =   \zeta\text{ on } \partial\Omega\label{eq2}
		\end{gather}
	\end{subequations}
	from the knowledge of $u$ in $\Omega=(0,1)\times(0,1)$, $\varphi\in L^2(\Omega)$ and $\zeta$ the trace of a function in $H^2(\Omega)$. This problem has been widely studied,  see for example \cite{rieder,kunisch}.

	The discretized version of the arising nonlinear least-squares problem is considered, obtained as described in \cite{rieder}. Namely problem  \eqref{eq1}-\eqref{eq2} was discretized using finite differences  choosing as grid points  $x_i=y_i=\frac{i-1}{N-1}$, for $i=1,\dots,N$ and $N=50$,    and  using lexicographical ordering, denoted by $l:\{1,\dots,N^2\}\rightarrow\{1,\dots,N^2\}$. Let   $\bar{\varphi}=[\bar \varphi_1,\dots,\bar\varphi_{N^2}]^T$,  where  $\bar\varphi_{l(i,j)}=\varphi(x_i,y_j)$ and 
	$A$ be the matrix arising from the discretization of the Laplacian operator. Moreover for $c\in\mathbb{R}^{N^2}$ let ${F}(c)=(A+\texttt{ diag}(c))^{-1}\bar{\varphi}$. Then, $F:\mathbb{R}^{N^2}\rightarrow\mathbb{R}^{N^2}$, and the resulting discrete problem is a  nonlinear least-squares problem  of size $n=N^2=2500$:  
	\begin{equation*}
	\min_{c\in\mathbb{R}^{N^2}} \frac{1}{2}\|F(c)-\bar u\|^2,
	\end{equation*}
	for a given $\bar{u}\in\mathbb{R}^{N^2}$. 
	For further details see \cite{rieder}.
	The tests were conducted choosing 
	$
	c(x,y)=1.5\sin(4\pi x)\sin(6\pi y) + 3 ((x-0.5)^2+(y-0.5)^2)+2
	$ as a parameter to be identified.
	The solution $u(x,y)$ of \eqref{eq} corresponding to this choice of $c(x,y)$ is 
	$
	u(x,y)=16 x(1-x)y(y-1)+1.
	$ Function 
	$\varphi$ in \eqref{eq} has been defined from \eqref{eq1}.   When the solution $u$ is analytically known, this is a zero residual problem. In order to obtain a nonzero residual problem the data $\bar{u}$ are artificially set as a perturbation of $[u_1,\dots,u_{N^2}]$ with $u_{l(i,j)}= u(x_i,y_j)$, to let $c^\dagger=[c^\dagger_1,\dots,c^\dagger_{N^2}]^T$, where $c^\dagger_{l(i,j)}= c(x_i,y_j)$, be a stationary point with strictly positive residual. Specifically 
	$\|J(c^\dagger)^T(F(c^\dagger)-\bar u)\|=0$ and $\|F(c^\dagger)-\bar u\|\simeq 0.1 $, for $J$ the Jacobian matrix of $F$.

	For this test problem the exact form of the the Jacobian matrix of $ F$ is given by:
	\begin{equation}\label{jacP2}
	J(c)=-(A+\texttt{diag}(c))^{-1}(\texttt{diag}( F(c))).
	\end{equation} 
	
	\item {\bf Problem 6.2} Given
	two images taken, for example, at different times, from different devices or perspectives, the goal is to determine a reasonable transformation, such that a transformed
	version of the first image is similar to the second one \cite{book_FAIR,haber2008adaptive,haber2006multilevel}. More specifically, given two images $T$ and $R$, the objective is to find a geometrical transformation $t$ (such as correction, deformation, displacement, distortion) such that the transformed image $T(t)$ is similar to $R$, where similarity is measured by a similarity/distance measure $D$. The transformation $t$ is recovered solving the following problem:
 \begin{equation*}
 \min_t D(T(t),R)+S(t),
  \end{equation*}
  with   $S$ a regularization term, which is usually added since 
 registration is an ill-posed problem. We will omit the regularization term and rely on the implicit regularization provided by our method. 
 Different choices are possible for $D$.  Here, the sum of squared differences (SSD) will be used \cite{book_FAIR}:
 \begin{equation*}
 D_{SSD}(T,R)=\frac{1}{2}\int_{\Omega}(T(x)-R(x))^2\,dx.
 \end{equation*}
 A discrete analogue of this distance, proposed in \cite[\S 6]{book_FAIR}, is given by a numerical approximation of the integral by a midpoint quadrature rule. Assume $\Omega\subset\mathbb{R}^d$ is divided into cells of width $h$ and cell centers $x_c$ and let $T_h=T(x_c)$, and $R_h = R(x_c)$, respectively. The discretized version of the SSD is defined as:
 \begin{equation*}
 \frac{1}{2}h\|T^h-R^h\|^2,\, h=\displaystyle\prod_{i=1}^{d}h_i.
 \end{equation*}
 To generate this test problem we used the code provided in the FAIR Matlab package\footnote{https://github.com/C4IR/FAIR.m}  and we obtain a nonlinear least-squares problem of $n=8320$ unknowns.

	\end{itemize} 
In the following, for uniformity of notation, for both tests it is assumed that the minimization problem to be solved is 
\begin{equation*}
\min_x\frac{1}{2}\|F(x)-y^\delta\|^2
\end{equation*}
and the true solution will be denoted by $x^\dagger$.

The practical implementation of the method will be now described.
All procedures were implemented in {\sc Matlab} and run using {\sc Matlab  2019a} on a MacBook Pro 2,4 GHz Intel Core i5, 4 GB RAM; the  machine precision is $\epsilon_m\sim 2\cdot 10^{-16}$. 
The trust-region procedure was implemented  according to
Algorithm \ref{algoTR}.

The major implementation issues are as follows. 
Lancsoz process in Algorithm \ref{algo_lanczos} is performed  to build  $Q_{\ell_k}$ and $T_{\ell_k}$ starting from $g_k$, and these matrices are used to compute the approximation to $p_k=B_k^{1/2}z_k$ in \eqref{plambda_inexact}. Experiments have shown that to get a good approximation of $p_k$, re-orthogonalization is necessary in the Lanczos process.  Re-orthogonalization also has a strong effect on the value of the quantity $\|Q_{\ell_k} T_{\ell_k}^TT_{\ell_k} Q_{\ell_k}^T-\ql\ql^TB_k\|$, which is required to be small to have a monotonic decrease of the error, cf. Lemma \ref{monotone_decay} and Remark \ref{remark_monotono}. We use the partial re-orthogonalization implemented in  \textit{Lanbpro}\footnote{https://github.com/epfl-lts2/unlocbox/blob/master/test$\_$bench/private/lanbpro.m} \cite{larsen}, \cite[\S 9.3.4]{gvl}. The resulting method will be labelled as LTR (Lanczos trust-region). 

The proposed method is going to be compared to its exact counterpart developed in \cite{ellittica}, that we  will label RTR (regularizing trust-region). For this, the square root of matrix $B_k$ is computed using the singular value decomposition of the Jacobian, provided by {\sc Matlab} function {\tt svd}.

Regarding the Jacobian matrix of $F$, the analytical expression was used for all test problems. Specifically, for {\bf Problem 6.1} the  Jacobian matrix has the form given  in \eqref{jacP2} and for {\bf Problem 6.2}  the analytic Jacobian is evaluated at each iterate  by the FAIR code.

In case of noisy problems, given the error level $\delta$, the exact data $y$ was perturbed by 
normally distributed values using
the {\sc Matlab} function {\tt randn}, in a way that $\|y-y^\delta\|=\delta$.

To compute the KKT point $(w_k,\lambda_k)$   at Step 4.1 \eqref{KKTy2} has to be solved. Since, from Lemma \ref{lemma_qcond}, the trust-region is ensured to be active,  this can be accomplished solving  the following  nonlinear scalar equation, called the \textit{secular equation}, bv Newton method  \cite[\S 7.3]{cgt}:
\begin{equation}\label{secular_eq}
\psi(\lambda)=0, \quad \psi(\lambda)= \frac{1}{\|w(\lambda)\|}-\frac{1}{\Delta_k}.
\end{equation} 
Starting from  an initial guess greater than the sought solution $\lambda_k$, the sequence generated converges monotonically to $\lambda_k$.  Typically, high accuracy in the solution of 
the above scalar equations is not needed, hence the Newton process  is terminated as soon as the
absolute value of function $\psi$ is below $10^{-2}$.
Each  Newton iteration applied to $\psi(\lambda)=0$  
requires the solution of a  linear system with shifted matrix of the form $(T_{\ell_k}^T T_{\ell_k})^2+\lambda I$.
The linear systems are solved  employing the singular value decomposition of $T_{\ell_k}$.

Algorithm \ref{algoTR} is run setting $\eta=10^{-1}$. In Step 3  the trust-region radius  is updated 
as follows 
\begin{eqnarray*}
	\Delta_0&=& \mu_0 \|\tilde{s}_0\| \qquad\qquad\;\;\;\mu_0=10^{-1}\\
	\Delta_{k+1}&=& \mu_{k+1} \|\tilde{s}_{k+1}\|   ,  \qquad  \  \mu_{k+1} = \begin{cases}  \displaystyle \frac{1}{6}\mu_k, & \text{ if } q_k<q \;\;\;\text{or} \;\;\; \rho_k<\eta_2, \\
		2\mu_k, & \text{ if } q_k>\nu q\;\; \mbox{and}\;\; \rho_k>\eta_2,\\
		\mu_k, & \text{otherwise}, \end{cases}
	\\%
\end{eqnarray*}
with $\tilde{s}_{k}$ defined in \eqref{tsk} and $q_k=\frac{\|Q_{\ell}^T(B_kp_k+g_k)\|}{\|g_k\|}$, $\nu=1.1$ and $\eta_2=0.25$, as in \cite{ellittica}.  We remark that $\|Q_{\ell}^T(B_kp_k+g_k)\|$  is computed as $\|T_{\ell}^T T_{\ell}Q_{\ell}^T p_k+\|g_k\|e_1\|$, with $e_1$ first vector of the canonical basis, so that  matrix $B_k$ does not need to be computed.  
The maximum and minimum values for $\Delta_k$ were set to
$\Delta_{\max}=10^4$ and  $\Delta_{\min}=10^{-12}$ and the maximum value for $\mu_k$ was set to $10^5$.
This updating strategy is inherited from \cite{ellittica} and is based on the following considerations. $\Delta_k$ given by the procedure described above 
preserves the property of converging to 
zero in case of exact data, as $\|\tilde{s}_k\|$ tends to zero. Further,
$\Delta_k$ is adjusted taking into account the  projected $q$-condition by monitoring  the value $q_k$, as it is satisfied whenever $q_k\geq q$. Therefore, if the  projected $q$-condition was not satisfied at the last computed iterate $\xkd$, 
it is reasonable to take a smaller
radius than in the case where it was fulfilled.
As was already observed in \cite{ellittica}, this updating strategy turns out to be efficient in practice. 

The free parameter $q$  was set equal to $0.8$, but the behaviour of the procedure  does not seem to be deeply affected by the value of $q$. 

The scalar $\bar\tau$ in the discrepancy principle \eqref{discrepance} cannot be chosen as required  in Lemma \ref{lemma_q_noise}, as in practice $\sigma$ is not known. We have chosen  $\bar{\tau}=0.1$ for {\bf Problem 6.1} and  $\bar{\tau}=10$ for {\bf Problem 6.2}.    We used two different values  in order to take into account that the unknown quantity $\sigma$ depends on the size of the problem, numerical tests provide an evidence of the effectiveness of this stopping rule. As $\tau_k$ depends on $k$, the stopping rule changes at each iteration. However, $\tau_k$ varies only slightly along the iterations as $\|J(x_k^\delta)\|$ is almost constant. Values of $\tau_k$ (computed at first iteration), are as follows:
\begin{itemize}
	\item {\bf Problem 6.1}:  $\tau_k\simeq 1.3e-3$,
	\item {\bf Problem 6.2}: $\tau_k\simeq 1.8e2$.
\end{itemize}

\subsection{Results on Problem 6.1}
In this section we focus on {\bf Problem 6.1}. We study the behaviour of the method depending on the choice of parameter $\ell_k$. In particular, we considered  constant values 
of $\ell_k$, i.e $\ell_k=\ell$, for each $k$, for different values of $\ell$, and an increasing choice of $\ell_k$ along the iterations. 
We measure the inexactness introduced by the use of Lanczos process in the approximation of the trust-region subproblem through \eqref{min_proj} and in the computation of $p_k$ via \eqref{def_p}.

We denote with $x^\dagger$ the exact solution of the problem, and we remind that $x^{\delta}_{k(\delta)}$ denotes the last computed iterate.

In  Tables \ref{tab1}-\ref{tab2} we report for $n=2500$ and $\delta=3.1e-2$:
  \begin{itemize}
  \item the relative error in the  computation of  $B_k^{1/2}g_k$: $err_r=\frac{\|\tilde{s}_k-B_k^{1/2}g_k\|}{\|B_k^{1/2}g_k\|}$ where $B_k^{1/2}g_k$ is computed by  the SVD of $J(\xkd)$ and $\tilde{s}_k$ is defined in \eqref{tsk};
  	\item the relative error $err_p=\frac{\|p_{k}-p_{ex}\|}{\|p_{ex}\|}$ between $p_{k}$, the inexact step computed at Step 4.2 of Algorithm \ref{algoTR}, and $p_{ex}$, the exact step arising from the solution of the trust-region problem \eqref{TR};
	\item {\sf it}, the number of nonlinear iterations;
	\item $\|F(x^{\delta}_{k(\delta)})-y^\delta\|$, the residual at the computed solution;
	\item {\sf it inner}, the number of iterations to solve the secular equation \eqref{secular_eq} for the minimization of the model;
	\item {\sf RMSE}=$\sqrt{\sum_{i=1}^{n}(x^\dagger(i)-x^{\delta}_{k(\delta)}(i))^2}$, root mean squared error between the computed solution and true solution $x^\dagger$;
	\item {\sf time(s)}, total CPU time (in seconds) required for the optimization process;
	\item  {\sf time(ratio)}, the ratio between CPU time required by RTR and CPU time required by LTR.
  \end{itemize}
  
    \begin{table}
	\centering
	\begin{tabular}{c|cccccc}
	\hline 
	&	 $\ell=5$ & $\ell=10$& $\ell=20$ & $\ell=40$ & $\ell=100$ & $\ell=3+\lceil k/2\rceil$  \\ 
	\hline
	$err_r=\frac{\|\tilde{s}_k-B_k^{1/2}g_k\|}{\|B_k^{1/2}g_k\|}$ 	&  6.e-5 & 7.e-6 & 6.e-7& 4.e-8 & 1.e-9& 8.e-5\\  
	$err_p=\frac{\|p_{k}-p_{ex}\|}{\|p_{ex}\|}$ 	& 7.e-14& 9.e-14& 3.e-13 & 2.e-13&2.e-13& 1.e-14\\  
	\end{tabular}
	\caption{{\bf Problem 6.1}, $n=2500$, $\delta=3.0e-2$, $p_{k}$ is the inexact step computed at Step 4.2 of Algorithm \ref{algoTR}, and $p_{ex}$ is the exact step arising from the solution of the trust-region problem \eqref{TR}.
}
	\label{tab1}
	\end{table}

   \begin{table}
	\centering
	\begin{tabular}{c|c|cccccc}
	\hline
	& RTR & &&& LTR &&\\
	\hline
	&  & $\ell=5$ & $\ell=10$ & $\ell=20$ & $\ell=40$ & $\ell=100$ & $\ell_k=3+\lceil k/2\rceil$ \\
	\hline
	\sf{it} &	 67&67&67&67&67&67&67 \\  
		$\|F(x^{\delta}_{k(\delta)})-y^\delta\|$ &	4.1e-2 &4.1e-2 &4.2e-2 &4.2e-2 &4.2e-2 &4.2e-2 & 4.1e-2  \\ 
		\sf{it inner} &	145&	147&	145&	145&	145&	145&	145 \\  
		\sf{err}	& 7.6e-1& 7.8e-1& 7.9e-1& 7.9e-1& 7.9e-1& 7.9e-1& 7.6e-1\\ 
		\sf{time(s)} & 2266 & 1138& 1142 & 1186 & 1225&1282 &1141\\
		\sf{time(ratio)} & 1 & 1.99& 1.98& 1.91 & 1.85&1.77 &1.99
			\end{tabular}
	\caption{{\bf Problem 6.1}, $n=2500$, $\delta=3.0e-2$.}
	\label{tab2}
\end{table}

\begin{figure}[htbp]
\begin{center}
\includegraphics[width=0.49\linewidth,height=0.40\linewidth]{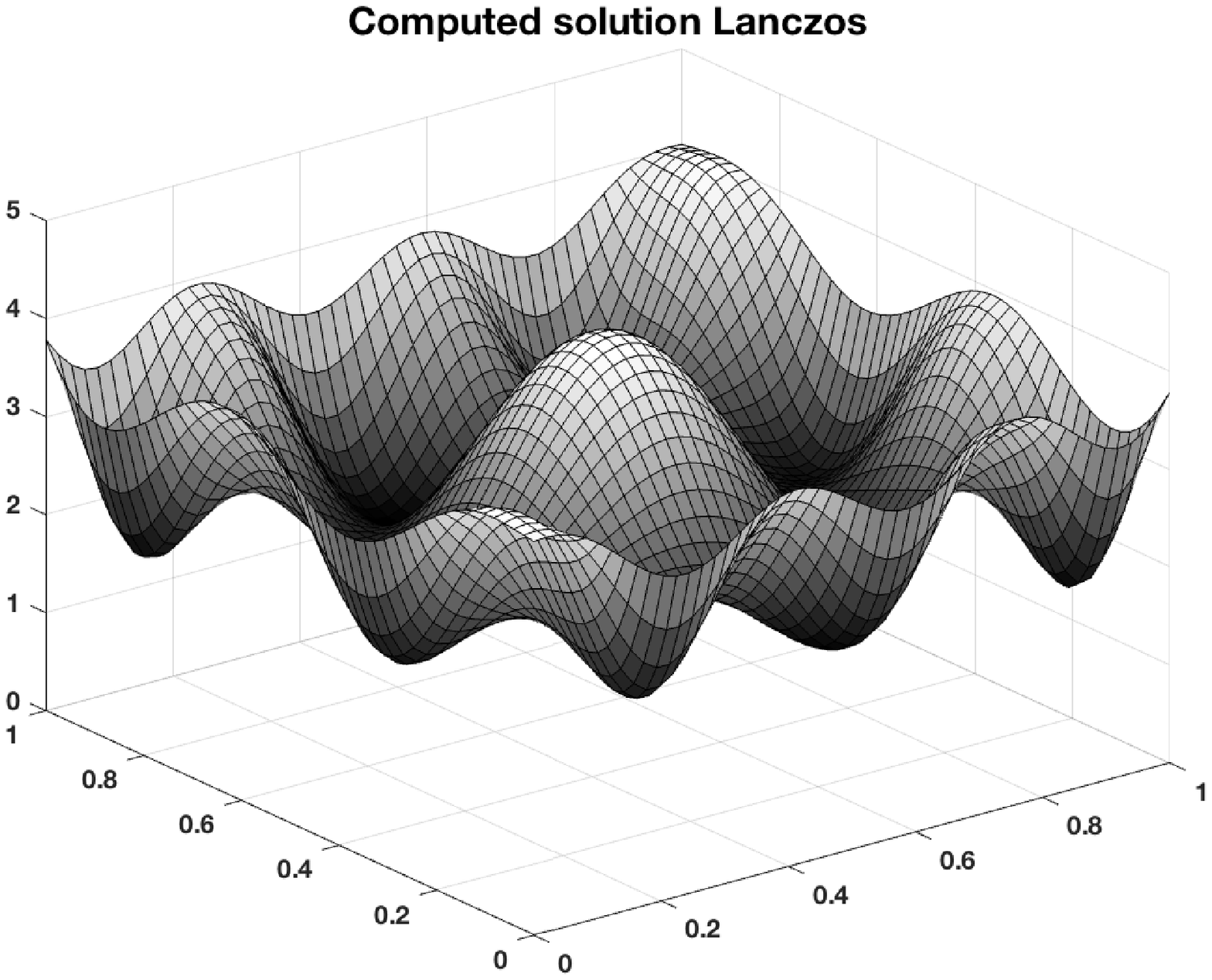}
\includegraphics[width=0.49\linewidth,height=0.40\linewidth]{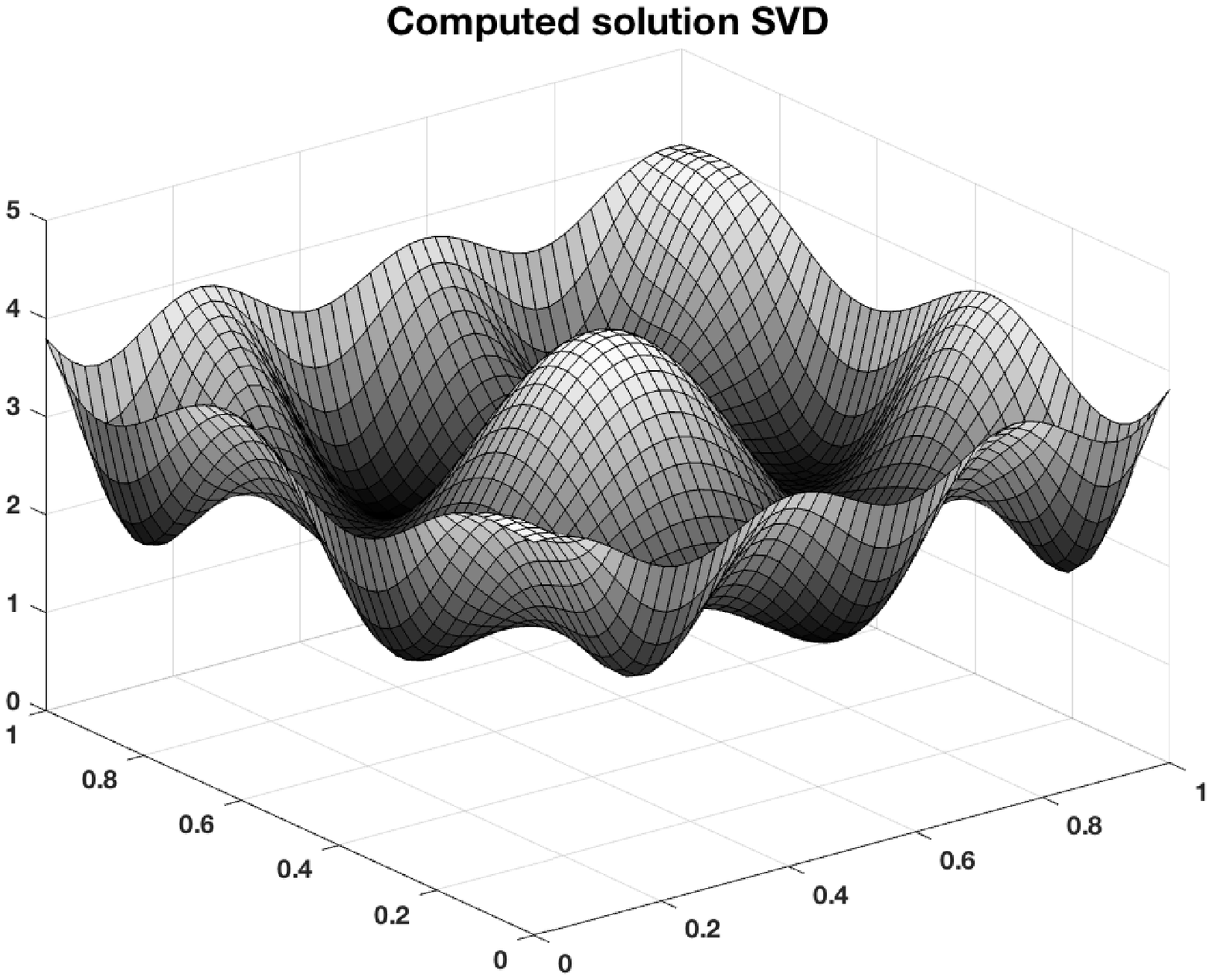}
\caption{{\bf Problem 6.1}, $n=2500$, $\delta=3.1e-2$. Plot of computed solution by LTR (left) and RTR (right).}
\label{fig:errsol}
\end{center}
\end{figure}

 From the Tables we remark that the value of $\ell$ affects the quality of the approximation of the right hand side, but it does not have an impact neither on the quality of the computed step nor  on the quality of the final solution approximation, as it is shown for example in Figure \ref{fig:errsol}. Moreover, the CPU times considerably decreases with $\ell$ and for all choices of $\ell$ it is significantly lower than the CPU time for the method employing the SVD decomposition. We considered also  the adaptive choice of $\ell_k=3+\lceil k/2\rceil$, which allows us to consider in the Lanczos process subspaces of increasing dimension, thus improving the precision toward the end of the optimization process, and which is more flexible, having the advantage of not requiring an a-priori choice for the parameter. 
 	
\begin{figure}[htbp]
\begin{center}
\includegraphics[width=0.49\linewidth,height=0.39\linewidth]{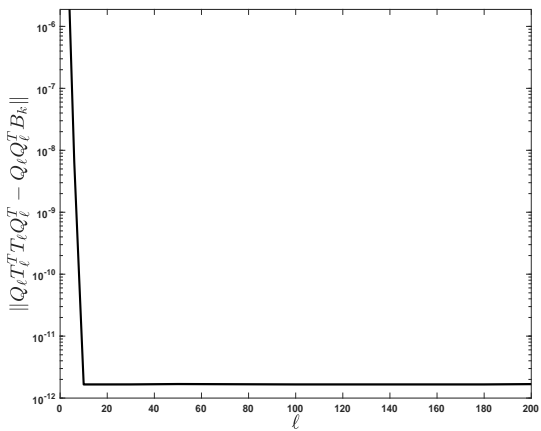}
\caption{{\bf Problem 6.1}, $n=2500$, $\delta=3.1e-2$. Plot of $\|Q_{\ell_k} T_{\ell_k}^TT_{\ell_k} Q_{\ell_k}^T-\ql\ql^TB_1\|$ as a function of $\ell$ .}
\label{fig:errsireortog}
\end{center}
\end{figure}

 In Figure \ref{fig:errsireortog} we report the plot of $\|Q_{\ell} T_{\ell}^TT_{\ell} Q_{\ell}^T-Q_{\ell}Q_{\ell}^TB_1\|$ as a function of $\ell$. The plot refers to  the first outer iteration of Algorithm \ref{algoTR}, but the behaviour is analogous also for the subsequent iterations.  
 We can notice that even for small  values of $\ell$ the quantity $\|Q_{\ell} T_{\ell}^TT_{\ell} Q_{\ell}^T-Q_{\ell}Q_{\ell}^TB_k\|$ is really small. Lemma \ref{monotone_decay} ensures that if this quantity is sufficiently small, the error at step $k$ decreases. 
 
\begin{figure}[htbp]
\begin{center}
\includegraphics[width=0.49\linewidth,height=0.39\linewidth]{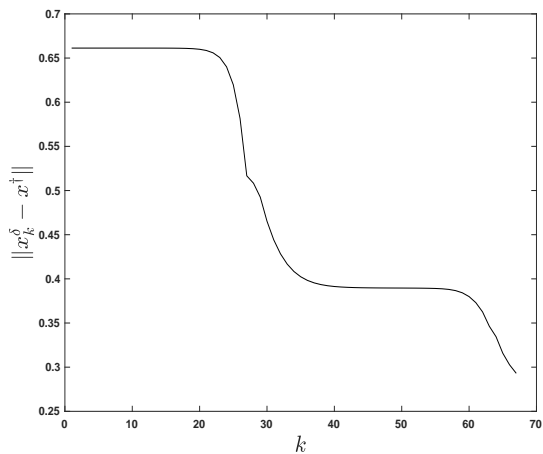}\\
\caption{{\bf Problem 6.1}, $n=2500$, $\delta=3.1e-2$. Plot of $\|\xkd-x^\dagger\|$ as a function of $k$ for $\ell_k=3+\lceil k/2\rceil$.}
\label{fig:errmonotono}
\end{center}
\end{figure}

 We show in Figure \ref{fig:errmonotono} that in practice the error is actually  monotonically decreasing. The figure refers to the adaptive  choice of $\ell_k$, $\ell_k=3+\lceil k/2\rceil$ but the behavior is the same for the other  choices of $\ell_k$. Then, in accordance with the theoretical results in the previous sections,  we recover the theoretical regularizing properties of RTR, interestingly even for small $\ell$ .

 From a practical point of view, there is then no evident benefit in using large values  of $\ell$ and an adaptive choice such as $\ell_k=3+\lceil k/2\rceil$ seems to be a good option.

\section{Results on Problem 6.2}
\begin{figure}[htbp]
\begin{center}
\includegraphics{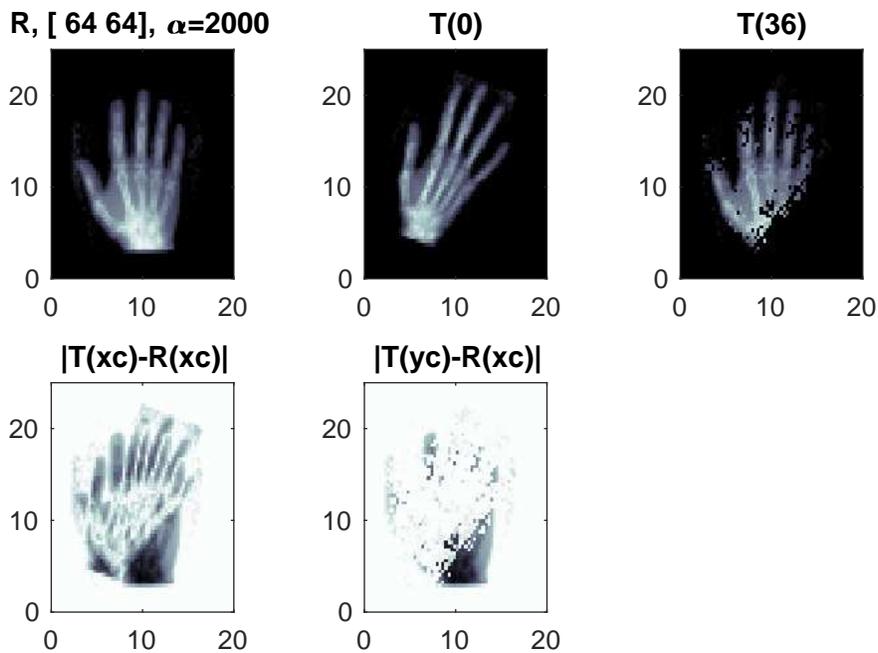}
\caption{{\bf Problem 6.2}, $n=8320$, $\delta=1.0e-1$. Plot of reconstructed image by method RTR.}
\label{fig:manosvd}
\end{center}
\end{figure}

In this section we report results on a test case taken from the FAIR library\footnote{E9$\_$Hands$\_$NPIRmb$\_$GN.m}. We consider two images of a hand, $R$ the reference image and $T$ the same image, rotated (up left and up center plots in Figure \ref{fig:manosvd} respectively).  The aim is to find a rotation to superimpose the two images. The domain is $\Omega=(0,20)\times(0,25)$ and spline interpolation, SSD distance  and elastic regularization are used, cf. \cite{book_FAIR}. The size of the problem is $n=8320$, and noise is added to image $T$ with $\delta=1.0e-1$.
The FAIR code implements a non-parametric image registration problem, solved by a Gauss-Newton (GN) method. We remark that the aim here is not to find the best method to solve the image registration problem, but rather to compare the performance of RTR and LTR on this example. 
 
\begin{figure}[htbp]
\begin{center}
\includegraphics{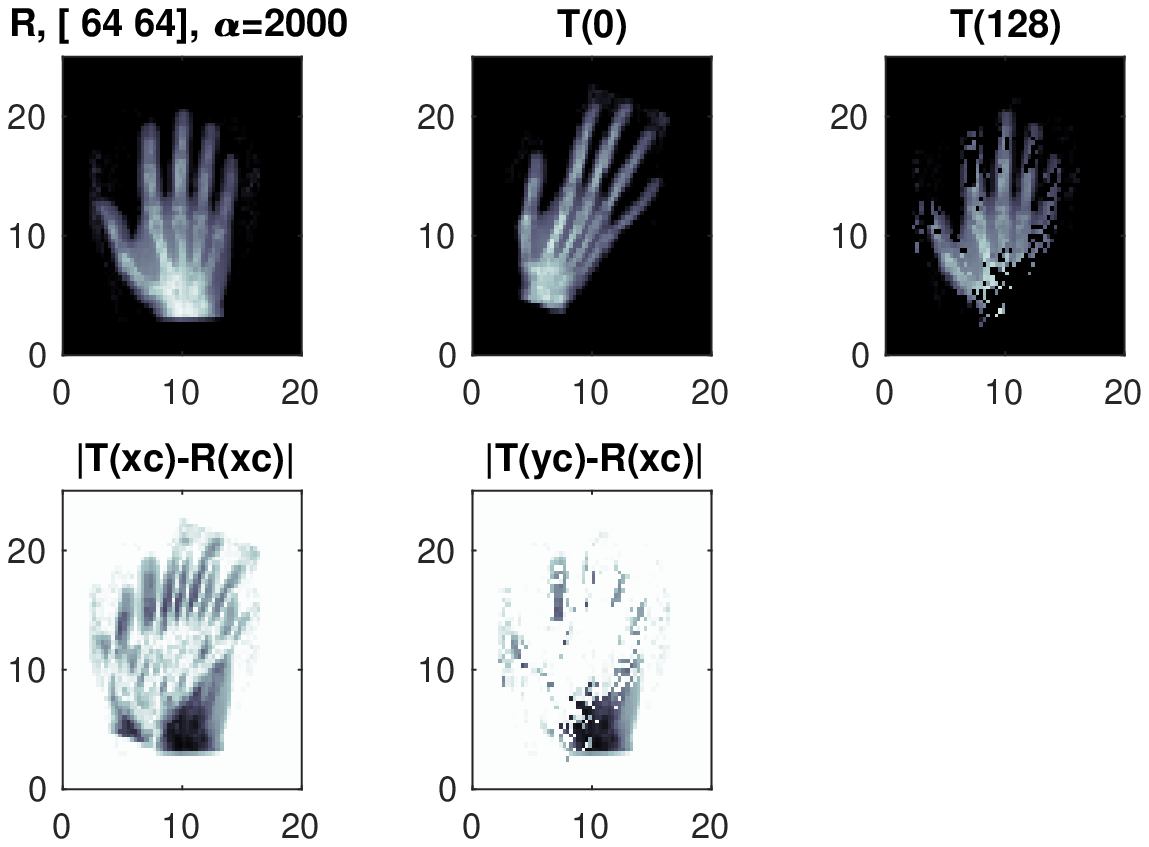}
\caption{{\bf Problem 6.2}, $n=8320$, $\delta=1.0e-1$.  Plot of reconstructed image by LTR with $\ell_k=100$ for all $k$.}
\label{fig:mano100}
\end{center}
\end{figure}

\begin{figure}[htbp]
\begin{center}
\includegraphics{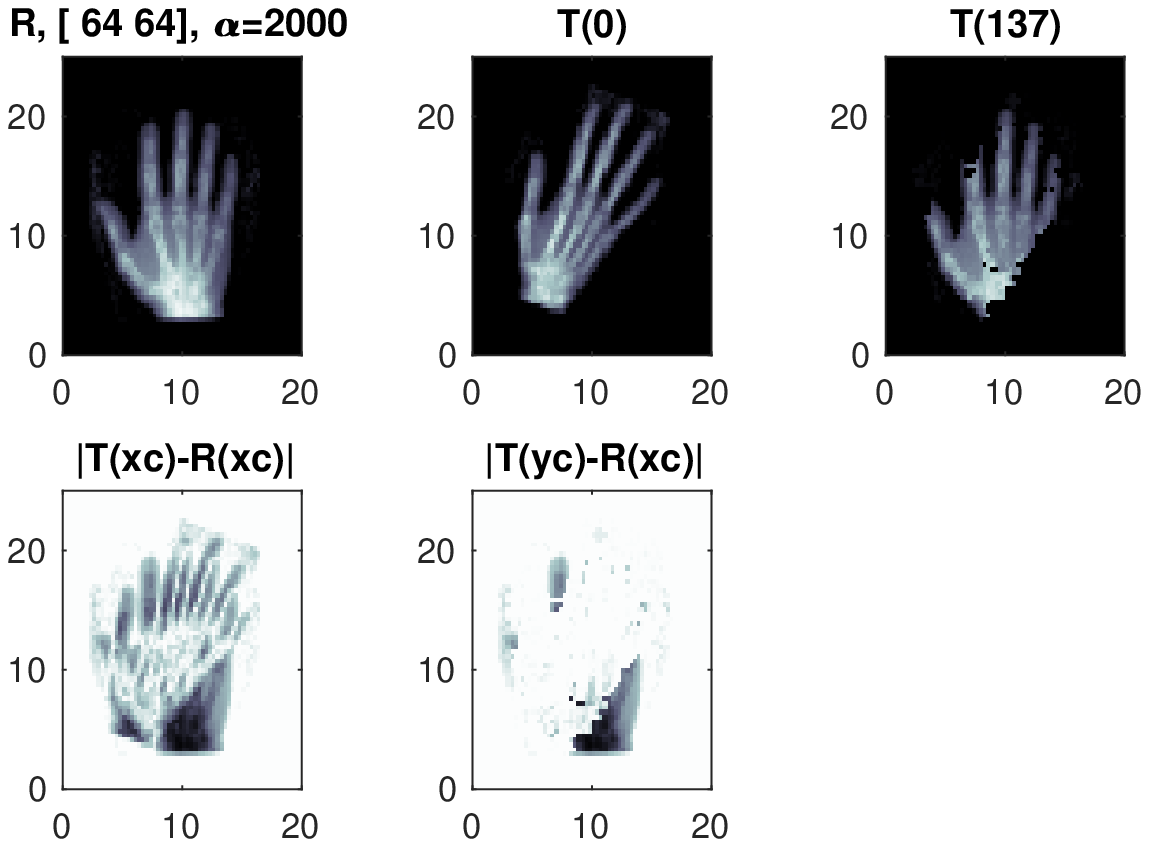}
\caption{{\bf Problem 6.2}, $n=8320$, $\delta=1.0e-1$.  Plot of reconstructed image by LTR with $\ell_k=10$ for all $k$.}
\label{fig:mano10}
\end{center}
\end{figure}

\begin{figure}[htbp]
\begin{center}
\includegraphics{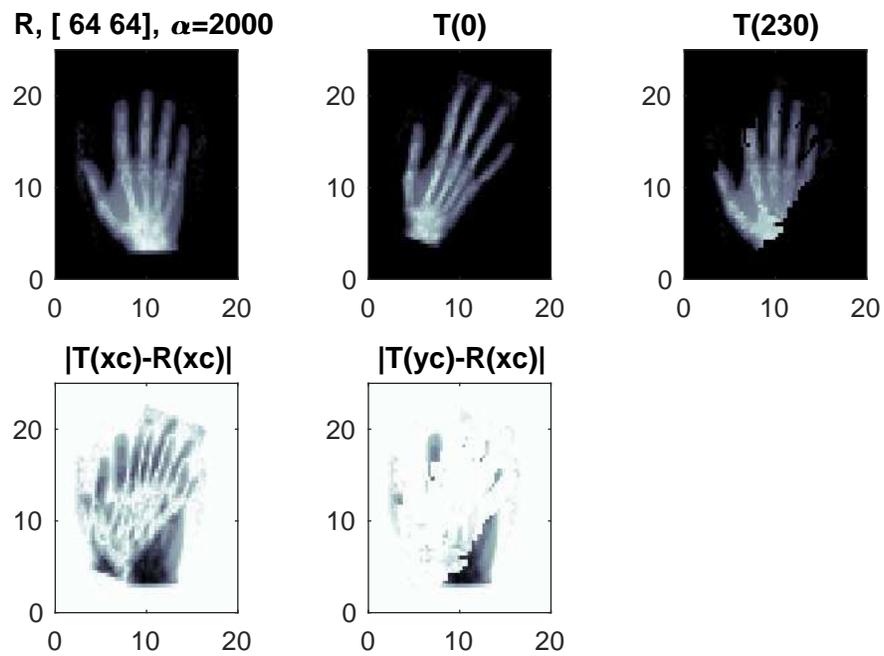}
\caption{{\bf Problem 6.2}, $n=8320$, $\delta=1.0e-1$.  Plot of reconstructed image by LTR with $\ell_k=3+\lceil k/2\rceil$.}
\label{fig:manok}
\end{center}
\end{figure}

\begin{figure}[htbp]
\begin{center}
\includegraphics{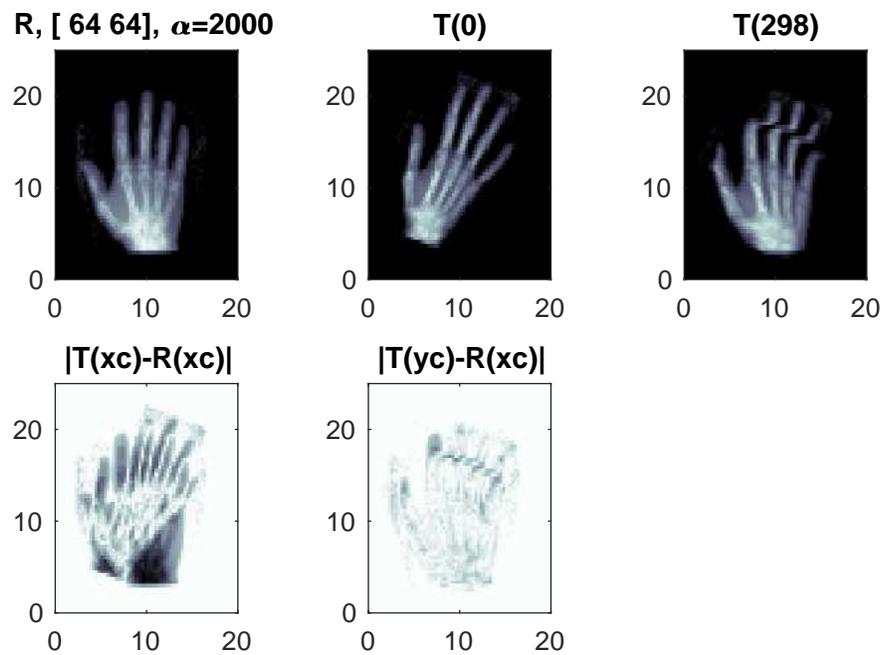}
\caption{{\bf Problem 6.2}, $n=8320$, $\delta=1.0e-1$. Plot of reconstructed image by Gauss-Newton method implemented in FAIR.}
\label{fig:manoGN}
\end{center}
\end{figure}

Image registration problems are usually really sensible to the choice of the regularization parameters, in our case the coefficient $\ell_k$. Common practice is to solve the problem for several values of the parameter and let specialists in the domain choose the best resulting image \cite{haber2008adaptive,haber2006multilevel}. In Figures \ref{fig:manosvd}-\ref{fig:manok} we report the reference image $R$ (up left), image $T$ at the beginning  (up center) and at the end (up right) of the optimization process, the difference between $R$ and $T$  at the beginning  (bottom left) and at the end (bottom center) of the optimization process. In Figure \ref{fig:manosvd} RTR is used, while in Figures \ref{fig:mano100} and  \ref{fig:mano10} we used  LTR with constant $\ell_k=\ell$ at each iteration, with $\ell=100,10$ respectively, and in Figure \ref{fig:manok} we used  LTR with the adaptive choice  $\ell_k=3+\lceil k/2\rceil$ as in the previous example. We can notice that the quality of the rotated image reconstructed using the SVD decomposition is not better than that of the image reconstructed by the proposed method. It is also interesting to notice that increasing $\ell$ this quality does not improve, actually lower values of $\ell$ provide the best results, as can be remarked also from the lower values achieved for the residual of the problem, as reported in Table \ref{tab4}. The use of a smaller projection space seems to improve the regularizing properties of the method. As for {\bf Problem 6.1}, an adaptive choice of $\ell$ seems to be a good choice to avoid hand-tuning. Also for this problem we can observe important time savings provided by the proposed Lanczos strategy. Iterations of LTR are so cheaper than those of RTR that, even if RTR takes far lesser  outer and inner iterations  to satisfy the discrepancy principle, the resulting computational time for LTR is much shorter. 
As a reference, we display in Figure \ref{fig:manoGN} the reconstructed image for the Gauss-Newton method originally proposed in the FAIR code. We can see that the lower part of the hand is better reconstructed, but the fingers are completely misinterpreted. 

For this problem the exact solution is not known, so we cannot verify if the error is monotonically decreasing. 

\begin{figure}[htbp]
\begin{center}
\includegraphics[scale=0.41]{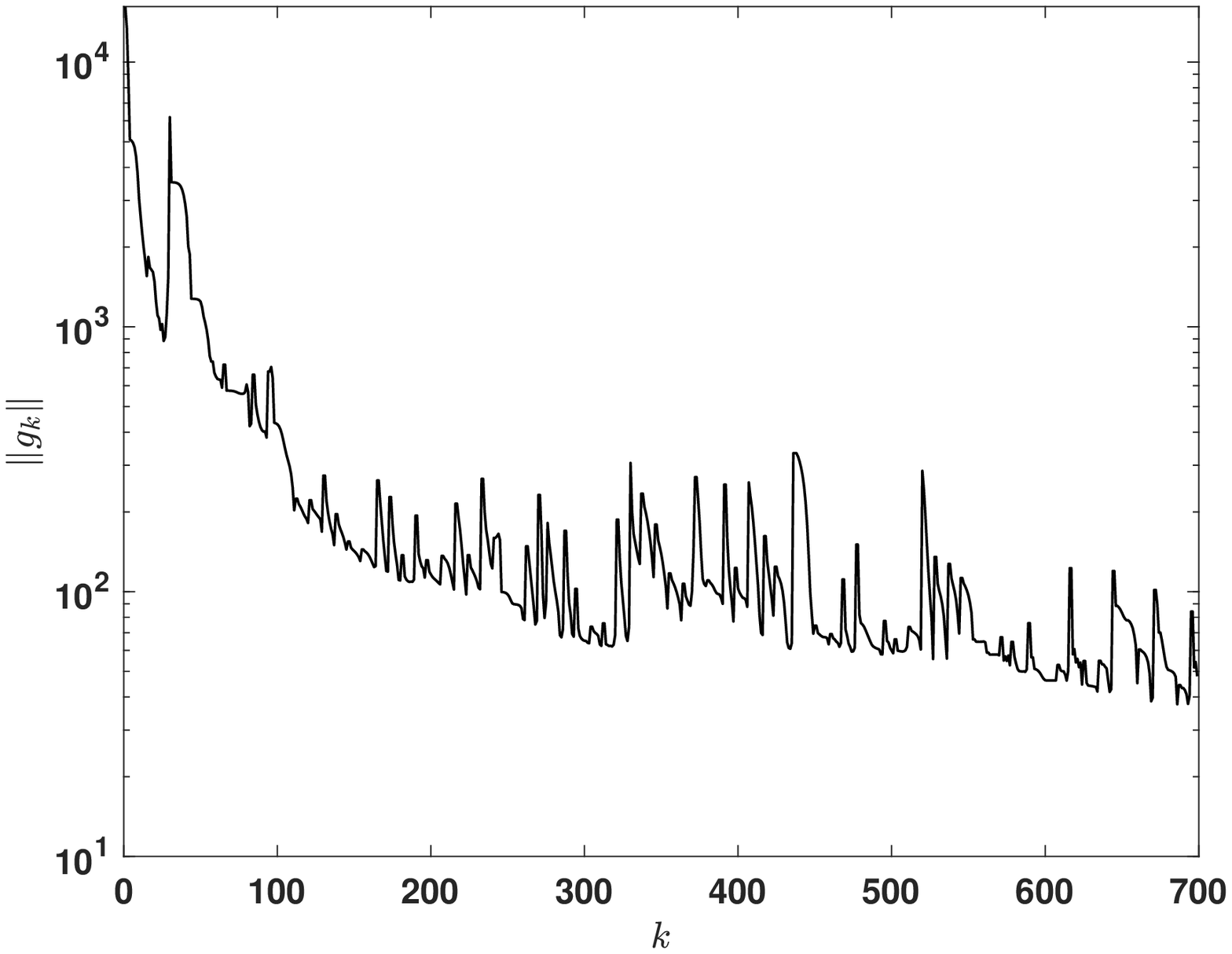}
\includegraphics[scale=0.41]{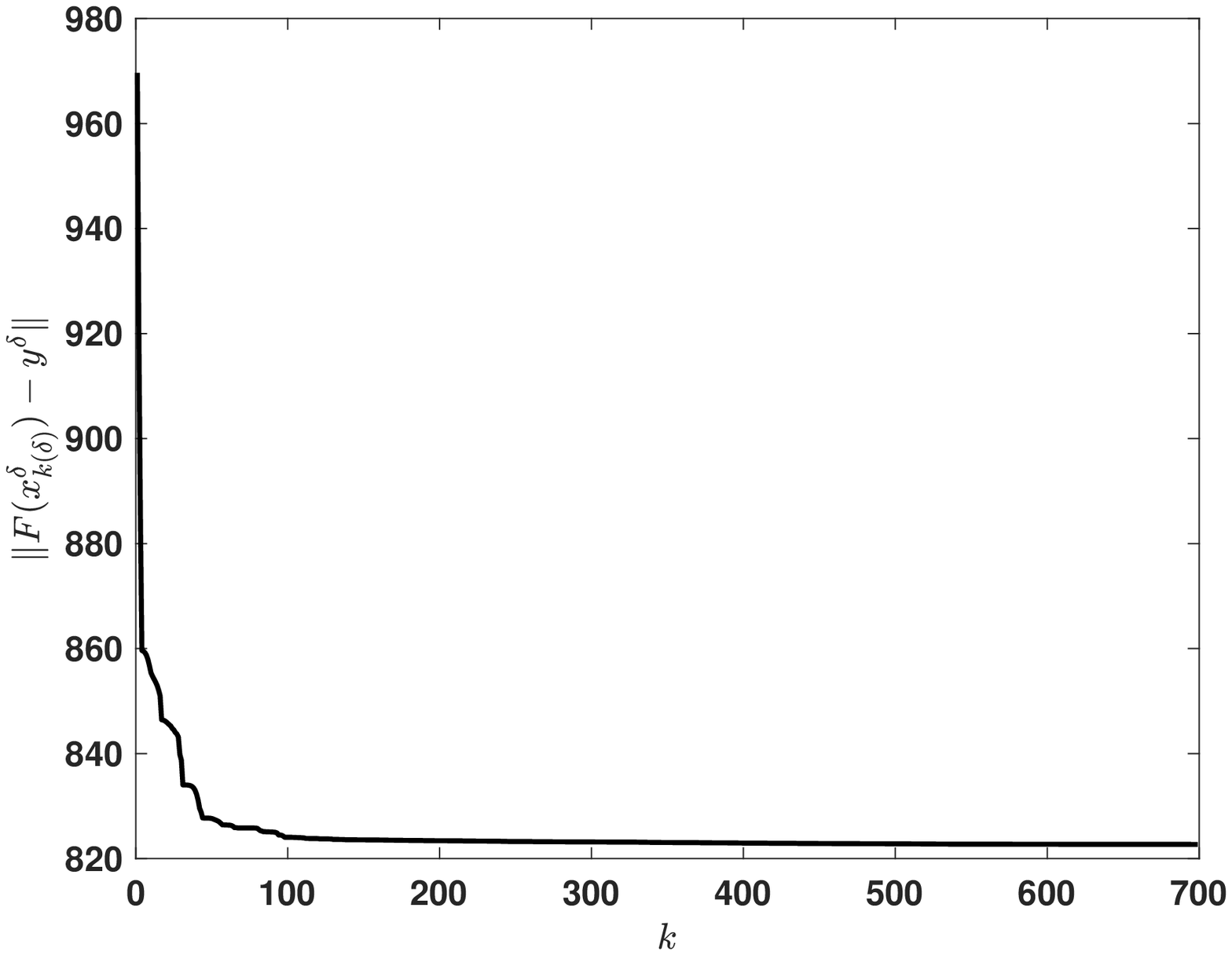}
\caption{{\bf Problem 6.2}, $n=8320$, $\delta=1.0e-1$. Plot of norm of the gradient (left) and norm of the residual (right) for LTR with $l=10$.}
\label{fig:norme}
\end{center}
\end{figure}

\begin{figure}[htbp]
\begin{center}
\includegraphics{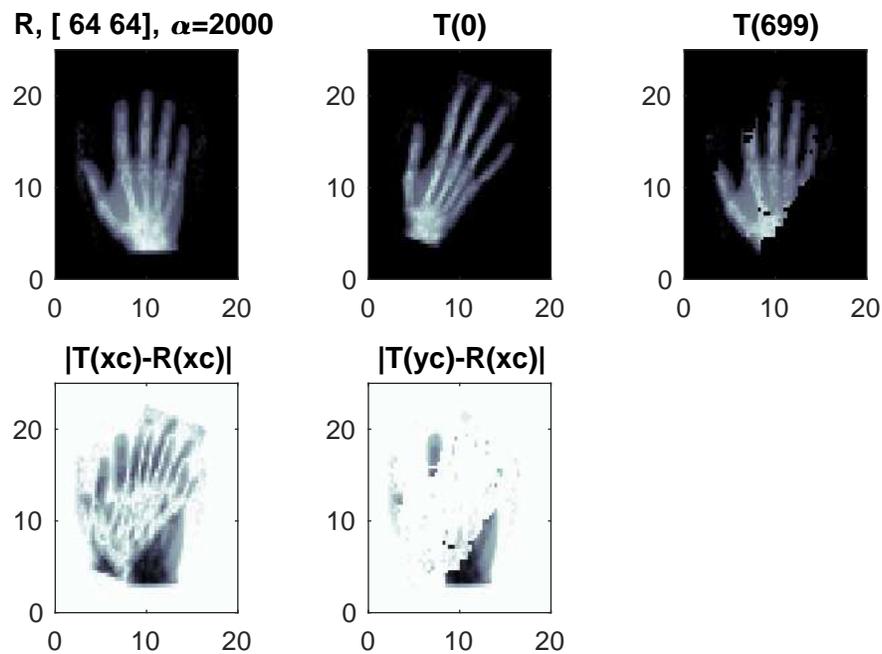}
\caption{{\bf Problem 6.2}, $n=8320$, $\delta=1.0e-1$. Plot of reconstructed image by LTR with $l=10$ and $\bar\tau=1$.}
\label{fig:mano_tau}
\end{center}
\end{figure}

In Figure \ref{fig:norme} we plot the norm of the gradient $\|g_k\|$ (left) and the norm of the residual $\|F(x^{\delta}_{k(\delta)})-y^\delta\|$ (right) along the iterations.  Just for this run we set $\bar\tau=1$. We can notice that the norm of the gradient is quite oscillatory and even allowing a higher number of iterations it does not decrease significantly and thus the stopping criterion \eqref{discrepance} is not met and  the method is stopped after the maximum number of iterations is reached. However, we can see that it is not useful to iterate further, as the norm of the objective function does not change and the quality of the solution approximation is not improved, cf. Figure \ref{fig:mano_tau}. Therefore, the value for $\bar\tau=10$ chosen for all the other runs seems to be a reasonable choice.

  \begin{table}
	\centering
	\begin{tabular}{c|c|cccc}
	\hline
	& RTR & & LTR &\\
	\hline
	& & $\ell=10$ & $\ell=50$ & $\ell=100$ & $\ell_k=3+\lceil k/2\rceil$ \\
	\hline
	\sf{it} &	 36&137&262&128&230\\  
		$\|F(x^{\delta}_{k(\delta)})-y^\delta\|$ &	9.1e2 &8.2e2 &9.3e2 &9.7e2 &7.7e2   \\ 
		\sf{it inner} &	162&	629&	1273&	684&	1123 \\  
		\sf{time(s)} & 10730 & 55& 214 & 243 & 273\\
		\sf{time(ratio)} & 1 & 195& 50 & 44 & 39
					\end{tabular}
	\caption{{\bf Problem 6.2}, $n=8320$, $\delta=1.0e-1$.}
	\label{tab4}
\end{table}

\section{Conclusions}
In this work we have proposed an extension of the method proposed in \cite{ellittica}, specially designed to handle large scale problems. The proposed approach is a hybrid Lanczos Tikhonov method, based on an elliptical trust-region implementation. An inexact solution of the trust-region subproblem is computed thanks to a Lanczos approach, which allows to considerably decrease the computational time without affecting the quality of the computed solution, as shown in the numerical tests. The complexity of the method is evaluated and 
some regularizing properties are proved in theory and verified in practice.

\bibliographystyle{plain}
\bibliography{biblio_ellittica_large}{}

\end{document}